# DYNAMICAL STABILITY OF PERCOLATION FOR SOME INTERACTING PARTICLE SYSTEMS AND $\varepsilon$-MOVABILITY

By Erik I. Broman[1] and Jeffrey E. Steif[2]

*Chalmers University of Technology*

In this paper we will investigate dynamic stability of percolation for the stochastic Ising model and the contact process. We also introduce the notion of downward and upward $\varepsilon$-movability which will be a key tool for our analysis.

**1. Introduction.** Consider bond percolation on an infinite connected locally finite graph $G$, where, for some $p \in [0,1]$, each edge (bond) of $G$ is, independently of all others, open with probability $p$ and closed with probability $1 - p$. Write $\pi_p$ for this product measure. The main questions in percolation theory (see [10]) deal with the possible existence of infinite connected components (clusters) in the random subgraph of $G$ consisting of all sites and all open edges. Write $\mathcal{C}$ for the event that there exists such an infinite cluster. By Kolmogorov's 0–1 law, the probability of $\mathcal{C}$ is, for fixed $G$ and $p$, either 0 or 1. Since $\pi_p(\mathcal{C})$ is nondecreasing in $p$, there exists a critical probability $p_c = p_c(G) \in [0,1]$ such that

$$\pi_p(\mathcal{C}) = \begin{cases} 0, & \text{for } p < p_c, \\ 1, & \text{for } p > p_c. \end{cases}$$

At $p = p_c$, we can have either $\pi_p(\mathcal{C}) = 0$ or $\pi_p(\mathcal{C}) = 1$, depending on $G$.

In [15] the authors initiated the study of dynamical percolation. In this model, with $p$ fixed, the edges of $G$ switch back and forth according to independent 2 state Markov chains where 0 switches to 1 at rate $p$ and 1 switches to 0 at rate $1 - p$. In this way, if we start with distribution $\pi_p$, the distribution of the system is at all times $\pi_p$. The general question studied in

Received June 2004; revised January 2005.

[1] Supported in part by the Swedish Natural Science Research Council.

[2] Supported in part by the Swedish Natural Science Research Council, NSF Grant DMS-01-0384 and the Göran Gustafsson Foundation (KVA).

[15] was whether there could exist atypical times at which the percolation structure looks different than at a fixed time.

We record here some of the results from [15]; (i) for any graph $G$ and for any $p < p_c(G)$, there are no times at which percolation occurs, (ii) for any graph $G$ and for any $p > p_c(G)$, there are no times at which percolation does not occur, (iii) there exist graphs which do not percolate for $p = p_c(G)$, but, nonetheless, for $p = p_c(G)$, there are exceptional times at which percolation occurs, (iv) there exist graphs which percolate for $p = p_c(G)$, but, nonetheless, for $p = p_c(G)$, there are exceptional times at which percolation does not occur, and (v) for $\mathbb{Z}^d$ with $d \geq 19$ with $p = p_c(\mathbb{Z}^d)$, there are no times at which percolation occurs. In addition, it has recently been shown in [23] that, for site percolation on the triangular lattice, for $p = p_c = 1/2$, there are exceptional times at which percolation occurs. Given this, a similar result would be expected for $\mathbb{Z}^2$.

The point of the present paper is to initiate a study of dynamical percolation for *interacting* systems where the edges or sites flip at rates which depend on the neighbors. We point out that in a different direction such questions in continuous space, but without interactions, related to continuum percolation have been studied in [2].

*Ising model results.* Precise definitions of the following Ising model measures and the stochastic Ising model will be given in Section 2. Fix an infinite graph $G = (S, E)$. Let $\mu^{+,\beta,h}$ be the plus state for the Ising model with inverse temperature $\beta$ and external field $h$ on $G$ [this is a probability measure on $\{-1, 1\}^S$]. Let $\Psi^{+,\beta,h}$ denote the corresponding stochastic Ising model; [this is a stationary continuous time Markov chain on $\{-1, 1\}^S$ with marginal distribution $\mu^{+,\beta,h}$]. Let $\mathcal{C}^+$ $(\mathcal{C}^-)$ denote the event that there exists an infinite cluster of sites with spin $1$ $(-1)$ and let $\mathcal{C}_t^+$ $(\mathcal{C}_t^-)$ denote the event that there exists an infinite cluster of sites with spin $1$ $(-1)$ at time $t$. It is known that the family $\mu^{+,\beta,h}$ is, for fixed $\beta$, stochastically increasing (to be defined later) in $h$.

THEOREM 1.1. *Consider a graph $G = (S, E)$ of bounded degree. Fix $\beta \geq 0$ and let $h_c = h_c(\beta)$ be defined by*

$$h_c := \inf\{h : \mu^{+,\beta,h}(\mathcal{C}^+) = 1\}.$$

*Then for all $h > h_c$,*

$$\Psi^{+,\beta,h}(\mathcal{C}_t^+ \text{ occurs for every } t) = 1$$

*and for all $h < h_c$,*

$$\Psi^{+,\beta,h}(\exists\, t \geq 0 : \mathcal{C}_t^+ \text{ occurs}) = 0.$$

*If we modify $h_c$ to be instead*

$$h'_c := \sup\{h : \mu^{+,\beta,h}(\mathcal{C}^-) = 1\},$$

*the same two claims hold with $\mathcal{C}_t^+$ replaced by $\mathcal{C}_t^-$ and with $h < h'_c$ and $h > h'_c$ reversed.*

This result tells us what happens in the subcritical and supercritical cases (with respect to $h$ with $\beta$ held fixed). It is the analogue of the easier Proposition 1.1 in [15] where it is proved that if $p < p_c$ ($p > p_c$), then, with probability 1, there is percolation at no time (at all times).

The following easy lemma gives us information about when $h_c$ is nontrivial.

LEMMA 1.2. *Assume the graph $G$ has bounded degree and let $\beta$ be arbitrary. Then $h_c > -\infty$. If $p_c(\text{site}) < 1$, then $h_c < \infty$. Similar results hold if $h_c$ is replaced by $h'_c$.*

The following theorems, where we restrict to $\mathbb{Z}^d$, will only discuss the case $h = 0$. However, this will in many cases give us information about the "critical" case $(\beta, h_c(\beta))$ since, in a number of situations, $h_c(\beta) = 0$. For example, this is true on all $\mathbb{Z}^d$ with $d \geq 2$ and $\beta$ sufficiently large. We also mention that while the relationship between $h_c$ and $h'_c$ in Theorem 1.1 might in general be complicated, for $\mathbb{Z}^d$, one easily has that $h_c = -h'_c$; this follows from the known fact that the plus and minus states are the same when $h \neq 0$. When $h = 0$, we will abbreviate $\mu^{+,\beta,0}$ by $\mu^{+,\beta}$ and $\Psi^{+,\beta,0}$ by $\Psi^{+,\beta}$. We point out that while $\mu^{+,\beta,h}$ is stochastically increasing in $h$ for fixed $\beta$, there is no such monotonicity in $\beta$ for fixed $h$, not even for $h = 0$. Therefore, we must use a different approach in the latter case.

We first study percolation of $-1$'s and then percolation of $1$'s. Let

$$\beta_p(2) := \inf\left\{\beta : \sum_{l=1}^{\infty} l3^{l-1}e^{-2\beta l} < \infty\right\} = \frac{\log 3}{2}.$$

We will refer to $\beta_p(2)$ as the critical inverse temperature of the Peierls regime for $\mathbb{Z}^2$. The choice of $\beta_p(2)$ might at first look quite arbitrary, but it is exactly what is needed to carry out a contour argument (known as Peierls argument) for $\mathbb{Z}^2$. For $d \geq 3$, there is a $\beta_p(d)$, such that, for $\beta$ larger than $\beta_p(d)$, a similar (although topologically more complicated) argument works for $\mathbb{Z}^d$. As a result of this "contour argument," it is well known and easy to show that, for $\beta > \beta_p(d)$, we have that

(1) $$\mu^{+,\beta}(\mathcal{C}^-) = 0.$$

Our next result is a dynamical version of (1) and we emphasize that this corresponds to the critical case as it is easy to check that, for these $\beta$'s, $h_c(\beta) = 0$.

THEOREM 1.3.  *For $\mathbb{Z}^d$ with $d \geq 2$ and $\beta > \beta_p(d)$,*

$$\Psi^{+,\beta}(\exists\, t \geq 0 : \mathcal{C}_t^-\ \ occurs) = 0.$$

It is well known that $\beta_p(d) \geq \beta_c(d)$, the latter being the critical inverse temperature for the Ising model on $\mathbb{Z}^d$. For $d = 2$, Theorem 1.3 can be extended down to the critical inverse temperature $\beta_c(2)$. First, it is known (see [5]) that on $\mathbb{Z}^2$, for *all* $\beta$,

$$(2) \qquad\qquad \mu^{+,\beta}(\mathcal{C}^-) = 0.$$

Our dynamical analogue for $\beta > \beta_c$ is the following where we again point out that this is also a critical case as it is easy to check that, for these $\beta$'s, we also have $h_c(\beta) = 0$.

THEOREM 1.4.  *For the stochastic Ising model $\Psi^{+,\beta}$ on $\mathbb{Z}^2$ with parameter $\beta > \beta_c$,*

$$\Psi^{+,\beta}(\exists\, t \geq 0 : \mathcal{C}_t^-\ \ occurs) = 0.$$

Interestingly, (1) is not always true for $\beta > \beta_c(d)$, although, as stated, it is true for $\mathbb{Z}^2$ or $\beta$ sufficiently large. In [1] it is shown that for $\mathbb{Z}^d$ with large $d$, there exists $\beta^+ > \beta_c(d)$ such that the probability in (1) is, in fact, 1 for all $\beta < \beta^+$. Moreover, they show that, for these $\beta$, there exists $h > 0$ with

$$\mu^{+,\beta,h}(\mathcal{C}^-) = 1.$$

For such $\beta$'s, this means that $h_c' > 0$ and, hence, it immediately follows from Theorem 1.1 that

$$\Psi^{+,\beta}(\mathcal{C}_t^-\ occurs\ for\ every\ t) = 1.$$

Note that, for these values of $\beta$, the case $h = 0$ is a noncritical case.

We next look at percolation of 1's under $\mu^{+,\beta}$. In the above results, we have not discussed the case of percolation of $-1$'s when $\beta \leq \beta_c$. However, by symmetry, this is the same as studying percolation of 1's in this case and so we can now move over to the study of $\mathcal{C}^+$.

First, it is well known that, for any graph of bounded degree, $\mu^{+,\beta,h} \neq \mu^{-,\beta,h}$ implies that $\mu^{+,\beta,h}(\mathcal{C}^+) = 1$. (This is proved in [3] for $\mathbb{Z}^d$; this argument extends to any graph of bounded degree.) In particular, for any graph $G$ of bounded degree and for $\beta > \beta_c(G)$,

$$(3) \qquad\qquad \mu^{+,\beta}(\mathcal{C}^+) = 1.$$

Our next result is a dynamical version of (3) for $\mathbb{Z}^d$. We mention that this result sometimes corresponds to a critical case and sometimes not. For $\beta > \beta_p(d)$ in $\mathbb{Z}^d$ or $\beta > \beta_c(2)$ in $\mathbb{Z}^2$, we have seen that $h_c = 0$ and so, in these cases, this next result covers the critical case. However, as pointed out, for $d$ large and $\beta$ just a little higher than $\beta_c$, the result in [1] gives us that $h_c < 0$ and, hence, in this case, this next theorem already follows from Theorem 1.1.

THEOREM 1.5. *For the stochastic Ising model $\Psi^{+,\beta}$ on $\mathbb{Z}^d$ with parameter $\beta > \beta_c(d)$,*

$$\Psi^{+,\beta}(\mathcal{C}_t^+ \text{ occurs for every } t) = 1.$$

(The proof we give actually works for any graph of bounded degree.) We mention that while $\beta > \beta_c$ is a sufficient condition for (3) to hold, it is certainly not necessary. For example, on $\mathbb{Z}^3$ we have that $\mu^{+,0}(\mathcal{C}^+) = 1$ since $\mu^{+,0} = \pi_{1/2}$ and the critical value for site percolation on $\mathbb{Z}^3$ is less than $1/2$. The reason $\beta_c$ appears is the connection between the Ising model and the random cluster model; $\beta_c$ corresponds to the critical value for percolation in the corresponding random cluster model (see [13]).

We are now left with the case $\beta \le \beta_c$. We will not be able to say too much since it is not known in all cases whether one has percolation at a fixed time. We first, however, have the following easy result for $d \ge 3$. We do not prove this result since it follows easily from the fact that the critical value for site percolation on $\mathbb{Z}^d$ is less than $1/2$ for $d \ge 3$, as this gives easily that $h_c(\beta) < 0$ for $\beta$ sufficiently small and, hence, Theorem 1.1 is applicable.

Note that the case $\beta = 0$ follows from the result in [15] mentioned above.

PROPOSITION 1.6. *For $d \ge 3$, there exists $\beta_1(d) > 0$ such that, for all $\beta < \beta_1(d)$, we have that*

$$\Psi^{+,\beta}(\mathcal{C}_t^+ \text{ occurs for every } t) = 1.$$

Finally, due to work of Higuchi, we can determine what happens with $\beta < \beta_c$ for $\mathbb{Z}^2$. It is shown in [16] that, for $\mathbb{Z}^2$, for all $\beta < \beta_c$, we have that $h_c(\beta) > 0$. The following result follows from this fact and Theorem 1.1.

THEOREM 1.7. *For $d = 2$, for all $\beta < \beta_c$, we have that*

$$\Psi^{+,\beta}(\exists t \ge 0 : \mathcal{C}_t^+ \text{ occurs}) = 0.$$

We note that even though it is known that for $\mathbb{Z}^2$, $\mu^{+,\beta_c}(\mathcal{C}^+) = 0$, we cannot conclude that

$$\Psi^{+,\beta_c}(\exists t \ge 0 : \mathcal{C}_t^+ \text{ occurs}) = 0,$$

since it is known (see [17]) that $h_c(\beta_c) = 0$. In contrast, based on the results in [23], it is interesting to ask the following:

QUESTION 1.8. *For the graph $\mathbb{Z}^2$, is it the case that*

$$\Psi^{+,\beta_c}(\exists t \geq 0 : \mathcal{C}_t^+ \text{ occurs}) = 1?$$

We finally mention that, interestingly, it is also known (see again [17]) that, for $\beta < \beta_c$, $\mu^{+,\beta,h_c(\beta)}(\mathcal{C}^+) = 0$.

*Contact process results.* Precise definitions of the following items will be given in Section 2. Fix an infinite graph $G = (S, E)$. Consider the contact process on $G = (S, E)$ with parameter $\lambda$. Denote by $\mu_\lambda$ the stochastically largest invariant measure, the so-called "upper invariant measure" (this is a probability measure on $\{0, 1\}^S$). Let $\Psi^\lambda$ denote the corresponding stationary contact process (this is a stationary continuous time Markov chain on $\{0, 1\}^S$ with marginal distribution $\mu_\lambda$). If $0 < \lambda_1 < \lambda_2$, it is well known that $\mu_{\lambda_1}$ is stochastically smaller than $\mu_{\lambda_2}$, denoted by

$$\mu_{\lambda_1} \preceq \mu_{\lambda_2}$$

(see Section 2 for this precise definition).

THEOREM 1.9. *Consider the contact process $\Psi^\lambda$ on a graph $G = (S, E)$, with initial and stationary distribution $\mu_\lambda$. Let $\lambda_p$ be defined by*

$$\lambda_p := \inf\{\lambda : \mu_\lambda(\mathcal{C}^+) = 1\}.$$

*We have that, for all $\lambda > \lambda_p$,*

$$\Psi^\lambda(\mathcal{C}_t^+ \text{ occurs for every } t) = 1.$$

In order for this theorem to be nonvacuous, we need to know that $\lambda_p < \infty$ for at least some graph. First, the fact that there exists $\lambda$ such that $\mu_\lambda(\mathcal{C}^+) > 0$ for $\mathbb{T}^d$ with $d \geq 2$ follows from [12]. Here $\mathbb{T}^d$ is the unique infinite connected graph without circuits and in which each site has exactly $d + 1$ neighbors; $\mathbb{T}^d$ is commonly known as the homogenous tree of order $d$. Combined with a 0–1 law which we develop, Proposition 4.2, we obtain that $\lambda_p < \infty$ in this case. For $\mathbb{Z}^d$ with $d \geq 2$ (as well as for $\mathbb{T}^d$), it is proved in [22] that, for large $\lambda$, $\mu_\lambda$ stochastically dominates high density product measures, which immediately implies that $\lambda_p < \infty$ in these cases.

When we prove Theorem 1.1, we will, in fact, prove a more general theorem which holds for a large class of systems. However, this proof will only work for models satisfying the so-called FKG lattice condition (which we call "monotone" in this paper). We now point out the important fact that, for $\lambda < 2$, in 1 dimension, the upper invariant measure for the contact process, while having positive correlations, is *not* monotone (see [20]). These terms are defined in Section 2. (One would also believe it is never monotone whenever the measure is not $\delta_0$.) Hence, Theorem 1.9 does not follow from the generalization of Theorem 1.1 which will come later.

*ε-movability.* We now introduce the concepts of upward and downward ε-movability. While we mainly introduce these as a technical tool to be used in our main results, it turns out that they are of interest in their own right. In [4] the concept of upward movability is studied for its own sake and is related to other well studied concepts, such as uniform insertion tolerance.

Let $S$ be a countable set. Take any probability measure $\mu$ on $\{-1,1\}^S$ and let $X$ be a $\{-1,1\}^S$ valued random variable with distribution $\mu$. Let $Z$ be a $\{-1,1\}^S$ valued random variable with distribution $\pi_{1-\varepsilon}$ and be independent of $X$. Define $X^{(-,\varepsilon)}$ by letting $X^{(-,\varepsilon)}(s) = \min(X(s), Z(s))$ for every $s \in S$, and let $\mu^{(-,\varepsilon)}$ denote the distribution of $X^{(-,\varepsilon)}$. In a similar way, define $X^{(+,\varepsilon)}$ by letting $X^{(+,\varepsilon)}(s) = \max(X(s), Z(s))$ for every $s \in S$, where $Z$ has distribution $\pi_\varepsilon$ and is independent of $X$. Denote the distribution of $X^{(+,\varepsilon)}$ by $\mu^{(+,\varepsilon)}$.

DEFINITION 1.10. Let $(\mu_1, \mu_2)$ be a pair of probability measures on $\{-1,1\}^S$, where $S$ is a countable set. Assume that

$$\mu_1 \preceq \mu_2.$$

If

$$\mu_1 \preceq \mu_2^{(-,\varepsilon)},$$

then we say that this pair of probability measures is downward ε-movable. If the pair is downward ε-movable for some $\varepsilon > 0$, we say that the pair is downward movable. Analogously, if

$$\mu_1^{(+,\varepsilon)} \preceq \mu_2,$$

then we say that the pair $(\mu_1, \mu_2)$ is upward ε-movable and that it is upward movable if the pair is upward ε-movable for some $\varepsilon > 0$.

For probability measures on $\{0,1\}^S$, we have identical definitions.

The relevance of downward (or upward) ε-movability to our dynamical percolation analysis will be explained in Section 5. In Section 3 we will prove ε-movability for general monotone systems, which will eventually lead to a proof of Theorem 1.1 (and its generalization). We now state a similar and key result for the contact process.

THEOREM 1.11. *Let $G$ be a graph of bounded degree, $0 < \lambda_1 < \lambda_2$ and $\mu_{\lambda_1}$, $\mu_{\lambda_2}$ be the upper invariant measures for the contact process on $\{0,1\}^S$ with parameters $\lambda_1$ and $\lambda_2$, respectively. Then $(\mu_{\lambda_1}, \mu_{\lambda_2})$ is downward movable.*

We finally mention how the above questions that we study fall into the context of classical Markov process theory. Let $(\Omega, \mathcal{F}, \mathbb{P})$ be the probability space where a stationary Markov process $\{X_t\}_{t \geq 0}$ taking values in some state space $\mathcal{S}$ is defined. Letting $\mu$ denote the distribution of $X_t$ (for any $t$), consider an event $\mathcal{A} \subseteq \mathcal{S}$ with $\mu(\mathcal{A}) = 1$. Let $\mathcal{A}_t$ be the event that $\mathcal{A}$ occurs at time $t$. We say that $\mathcal{A}$ is a *dynamically stable* event if $\mathbb{P}(\mathcal{A}_t \ \forall t \geq 0) = 1$. In Markov process terminology, this is equivalent to saying that $\mathcal{A}^c$ has *capacity zero*. All the questions in this paper deal with showing, for various models and parameters, that the event that there exists/there does not exist an infinite connected component of sites which are all open is dynamically stable.

The rest of this paper is divided into 9 sections. In Section 2 we will give all necessary preliminaries and precise definitions of our models. Sections 3 and 4 will deal with the concept of $\varepsilon$-movability. In Section 3 we develop what will be needed to prove Theorem 1.1 and its generalization. In Section 4 we will prove Theorem 1.11 (which is the key to Theorem 1.9), as well as give a proof that $\lambda_p < \infty$ for trees. In Section 5 we prove two elementary lemmas which relate the notion of $\varepsilon$-movability to dynamical questions. In the remaining sections proofs of the remaining results are given. We note that the proof of Theorem 1.4 will use the proof of Theorem 1.5 and, hence, will come afterward.

We end with one bit of notation. If $\mu$ is a probability measure on some set $U$, we write $X \sim \mu$ to mean that $X$ is a random variable taking values in $U$ with distribution $\mu$.

## 2. Models and definitions.
Before presenting the interacting particle systems discussed in this paper, we will present some definitions and results related to stochastic domination. Let $S$ be any countable set. For $\sigma, \sigma' \in \{-1, 1\}^S$, we write $\sigma \preceq \sigma'$ if $\sigma(s) \leq \sigma'(s)$ for every $s \in S$. An increasing function $f$ is a function $f : \{-1, 1\}^S \to \mathbb{R}$ such that $f(\sigma) \leq f(\sigma')$ for all $\sigma \preceq \sigma'$. For two probability measures $\mu, \mu'$ on $\{-1, 1\}^S$, we write $\mu \preceq \mu'$ if, for every continuous increasing function $f$, we have that $\mu(f) \leq \mu'(f)$. [$\mu(f)$ is shorthand for $\int f(x) \, d\mu(x)$.] When $\{-1, 1\}^S$ is replaced by $\{0, 1\}^S$, we have identical definitions. Strassen's theorem (see [19], page 72) states that if $\mu \preceq \mu'$, then there exist random variables $X, X'$ with distribution $\mu, \mu'$, respectively, such that $X \preceq X'$ a.s.

A very useful result is the so-called Holley's inequality, which appeared first in [18]. We will present a variant of the theorem by Holley; it is not the most general, but is sufficient for our purposes.

THEOREM 2.1. *Take $S$ to be a finite set. Let $\mu$, $\mu'$ be probability measures on $\{-1, 1\}^S$ which assign positive probability to all configurations $\sigma \in$*

$\{-1, 1\}^S$. *Assume that*

$$\mu(\sigma(s) = 1 | \sigma(S \setminus s) = \xi) \le \mu'(\sigma(s) = 1 | \sigma(S \setminus s) = \xi')$$

*for every $s \in S$ and $\xi \preceq \xi'$, where $\xi, \xi' \in \{-1, 1\}^{S \setminus s}$. Then $\mu \preceq \mu'$.*

PROOF. See [9] or [13] for a proof. $\square$

Two properties of probability measures which are often encountered within the field of interacting particle systems are the monotonicity property and the property of positive correlations presented below.

DEFINITION 2.2. Take $S$ to be a finite set. A probability measure $\mu$ on $\{-1, 1\}^S$ which assigns positive probability to every $\sigma \in \{-1, 1\}^S$ is called monotone if, for every $s \in S$ and $\xi \preceq \xi'$ where $\xi, \xi' \in \{-1, 1\}^{S \setminus s}$,

$$\mu(\sigma(s) = 1 | \sigma(S \setminus s) = \xi) \le \mu(\sigma(s) = 1 | \sigma(S \setminus s) = \xi').$$

We point out immediately that it is known that this is equivalent to the so-called FKG lattice condition.

DEFINITION 2.3. A probability measure $\mu$ on $\{-1, 1\}^S$ is said to have positive correlations if, for all bounded increasing functions $f, g: \{-1, 1\}^S \to \mathbb{R}$, we have

$$\mu(fg) \ge \mu(f)\mu(g).$$

The following important result is sometimes known as the FKG inequality (see [7]).

THEOREM 2.4. *Take $S$ to be a finite set. Let $\mu$ be a monotone probability measure on $\{-1, 1\}^S$ which assigns positive probability to every configuration. Then $\mu$ has positive correlations.*

PROOF. This was originally proved in [7]; see also [9] for a proof. $\square$

In this section, and also later in this paper, we will talk about convergence of probability measures. Convergence will always mean weak convergence, where $\{0, 1\}^S$ is given the product topology.

2.1. *The Ising model.* Take $G = (S, E)$, where $|S| < \infty$. The Ising measure $\mu^{\beta,h}$ on $\{-1, 1\}^S$ at inverse temperature $\beta \geq 0$, external field $h$ and with free boundary conditions is defined as follows. For any configuration $\sigma \in \{-1, 1\}^S$, let

$$(4) \qquad H^{\beta,h}(\sigma) = -\beta \sum_{\substack{\{t,t'\} \in E \\ t,t' \in S}} \sigma(t)\sigma(t') - h \sum_{t \in S} \sigma(t).$$

$H^{\beta,h}$ is called the Hamiltonian. Define $\mu^{\beta,h}$ by assigning the probability

$$(5) \qquad \mu^{\beta,h}(\sigma) = \frac{e^{-H^{\beta,h}(\sigma)}}{Z}$$

to any configuration $\sigma \in \{-1, 1\}^S$, where $Z$ is a normalization constant. Of course, $Z$ depends on the graph and the values $\beta$ and $h$, but this will not be important for us and, therefore, not reflected in the notation.

Take $S_n = \Lambda_{n+1} = \{-n-1, \ldots, n+1\}^d$ and $E_n$ to be the set of all nearest neighbor pairs of $S_n$. Given a configuration $\xi$ on $\{-1, 1\}^{\mathbb{Z}^d \setminus \Lambda_n}$, let, for $\sigma \in \{-1, 1\}^{\Lambda_n}$,

$$(6) \quad H_n^{\xi,\beta,h}(\sigma) = -\beta \sum_{\substack{\{t,t'\} \in E_n \\ t,t' \in \Lambda_n}} \sigma(t)\sigma(t') - h \sum_{t \in \Lambda_n} \sigma(t) - \beta \sum_{\substack{\{t,t'\} \in E_n \\ t \in \Lambda_n \\ t' \in \Lambda_{n+1} \setminus \Lambda_n}} \sigma(t)\xi(t')$$

be our Hamiltonian. Here $\xi$ is called a boundary condition. Again, we define a probability measure using (5), but using the Hamiltonian of (6) instead. This Ising measure will be denoted by $\mu_n^{\xi,\beta,h}$. The cases $\xi \equiv 1$ and $\xi \equiv -1$ are especially important and the corresponding Ising measures are denoted by $\mu_n^{+,\beta,h}$ and $\mu_n^{-,\beta,h}$, respectively. We view $\mu_n^{+,\beta,h}$ (resp. $\mu_n^{-,\beta,h}$) as a probability measure on $\{-1, 1\}^{\mathbb{Z}^d}$ by letting, with probability 1, the configuration be identically 1 (resp. $-1$) outside $\Lambda_n$. It is known (see [19], page 189) that the sequences $\{\mu_n^{+,\beta,h}\}$ and $\{\mu_n^{-,\beta,h}\}$ converge as $n$ tends to infinity; these limits are denoted by $\mu^{+,\beta,h}$ and $\mu^{-,\beta,h}$.

The same kind of construction can be carried out on any infinite connected locally finite graph $G = (S, E)$. One defines a Hamiltonian analogous to the one in (6), but with $\Lambda_n$ replaced by any $\Lambda \subseteq S$ where $|\Lambda| < \infty$. With $\xi \equiv 1$ or $\xi \equiv -1$, one then considers the corresponding limits of Ising measures as $\Lambda \uparrow S$, the limit turning out to be independent of the particular choice of sequence. See, for instance, [9] for how this is carried out in detail. Fix $h = 0$ and abbreviate $\mu^{+,\beta,0}$ and $\mu^{-,\beta,0}$ by $\mu^{+,\beta}$ and $\mu^{-,\beta}$. It is well known [8, 9] that, for any graph, there exists $\beta_c \in [0, \infty]$ such that, for $0 \leq \beta < \beta_c$, we have that $\mu^{-,\beta} = \mu^{+,\beta}$ (and there is then a unique so-called Gibbs state) and for $\beta > \beta_c$, $\mu^{-,\beta} \neq \mu^{+,\beta}$. For $\mathbb{Z}^d$ with $d \geq 2$, and many other graphs, $\beta_c \in (0, \infty)$. $\beta_c$ is sometimes referred to as the critical inverse temperature for phase transition in the Ising model. Furthermore, in [14] the author shows that if $G$ is of bounded degree, the condition $\beta_c < \infty$ is equivalent to the condition $p_c < 1$, where $p_c$ is the critical parameter value for site percolation on $G$. It is easy to see that for any graph of bounded degree $p_c > 0$ (see the proof of Theorem 1.10 of [10]). This, in turn, implies, via the connection between the random cluster model and the Ising model described below, that $\beta_c > 0$ for any graph of bounded degree.

2.2. *Spin systems.* A configuration $\sigma \in \{-1, 1\}^S$ can be seen as particles on a discrete set $S$ having one of two different "spins" represented by $-1$ and $1$. To this we will add a stochastic dynamics, and assume that the system is described by "flip rate intensities," which we will denote by $\{C(s, \sigma)\}_{s \in S,\ \sigma \in \{-1,1\}^S}$. $C(s, \sigma)$ represents the rate at which site $s$ changes its state when the present configuration is $\sigma$. Of course, $C(s, \sigma) \geq 0\ \forall s \in S, \sigma \in \{-1, 1\}^S$, and we assume that the interaction is nearest neighbor in the sense that the flip rate of a site $s \in S$ only depends on the configuration $\sigma$ at $s$ and at sites $t$ with $\{s, t\} \in E$. We will limit ourselves to only allow one site flip in every transition and we will only consider flip rate intensities such that

$$\sup_{s, \sigma} C(s, \sigma) < \infty.$$

In many cases we will consider translation invariant systems and then this last condition will hold trivially. Furthermore, we will always assume the trivial condition that, for every $s \in S$,

$$\sup_{\sigma\, :\, \sigma(s)=0} C(s, \sigma(s)) > 0, \qquad \sup_{\sigma\, :\, \sigma(s)=1} C(s, \sigma(s)) > 0.$$

We will call such an object a spin system (see [6] or [19] for results concerning general spin systems). Given such rates, one can obtain a Markov process $\Psi$ on $\{-1, 1\}^S$ governed by these flip rates; see [19]. Such a Markov process with a specified initial distribution $\mu$ on $\{-1, 1\}^S$ will be denoted by $\Psi^\mu$. Given a Markov process, $\mu$ will be called an invariant distribution for the process if the projections of $\Psi^\mu$ onto $\{-1, 1\}^S$ at any fixed time $t \geq 0$ is $\mu$. In this case, $\Psi^\mu$ will be a stationary Markov process on $\{-1, 1\}^S$, all of whose marginal distributions are $\mu$. Of course, the state space $\{-1, 1\}^S$ can be exchanged for either $\{0, 1\}^S$ or $\{0, 1\}^E$.

Sometimes we will work with two different sets of flip rates, $\{C_1(s, \sigma)\}_{s \in S, \sigma \in \{-1,1\}^S}$ and $\{C_2(s, \sigma)\}_{s \in S, \sigma \in \{-1,1\}^S}$, governing two Markov processes $\Psi_1$ and $\Psi_2$, respectively. We will write $C_1 \preceq C_2$ if the following conditions are satisfied:

(7) $\quad C_2(s, \sigma_2) \geq C_1(s, \sigma_1) \qquad \forall\, s \in S,\ \forall\, \sigma_1 \preceq \sigma_2 \text{ s.t. } \sigma_1(s) = \sigma_2(s) = 0,$

and

$$(8) \quad C_1(s, \sigma_1) \geq C_2(s, \sigma_2) \qquad \forall s \in S, \ \forall \sigma_1 \preceq \sigma_2 \text{ s.t. } \sigma_1(s) = \sigma_2(s) = 1.$$

The point of $C_1 \preceq C_2$ is that a coupling of $\Psi_1$ and $\Psi_2$ will then exist for which $\{(\eta, \delta) : \eta(s) \leq \delta(s) \ \forall s \in S\}$ is invariant for the process; see [19].

2.3. *Stochastic Ising models.* We will now briefly discuss stochastic Ising models. We will omit most details; for an extensive discussion and analysis, see again [19]. Consider $G_n = (S_n, E_n)$, defined in Section 2.1. Given $\beta$ and $h$, it is possible to construct flip rates $C_n^+$ on $\{-1, 1\}^{S_n}$ for which $\mu_n^{+,\beta,h}$ is reversible and invariant. We denote by $\Psi_n^{+,\beta,h}$ the corresponding stationary Markov process with initial distribution $\mu_n^{+,\beta,h}$. One possible choice of flip rate intensities are that, for every $s \in \Lambda_n$ and $\sigma \in \{-1, 1\}^S$,

$$\begin{aligned} &C_n^+(s, \sigma) \\ &= \exp\left[-\beta\left(\sum_{t \in \Lambda_n : \{t,s\} \in E_n} \sigma(t)\sigma(s) + \sum_{t \in \Lambda_{n+1} \setminus \Lambda_n : \{t,s\} \in E_n} \sigma(s)\right) - h\sigma(s)\right]. \end{aligned}$$

Sites in $\Lambda_{n+1} \setminus \Lambda_n$ are kept fixed at 1. Observe that if $s \in \Lambda_{n-1}$, the second sum is over an empty set. A straightforward calculation gives

$$(9) \qquad C_n^+(s, \sigma)\mu_n^{+,\beta,h}(\sigma) = C_n^+(s, \sigma_s)\mu_n^{+,\beta,h}(\sigma_s),$$

where

$$\sigma_s(t) = \begin{cases} \sigma(t), & \text{if } t \neq s, \\ -\sigma(t), & \text{if } t = s. \end{cases}$$

This shows that indeed $\mu_n^{+,\beta,h}$ is reversible and invariant for $C_n^+$. Any family of spin rates satisfying (9) is called a stochastic Ising model (on our finite set). One can show that there exists a limiting distribution $\Psi^{+,\beta,h}$ of $\Psi_n^{+,\beta,h}$ when $n$ tends to infinity; see [19], Theorem 2.2, page 17 and Theorem 2.7, page 139. Furthermore, $\Psi^{+,\beta,h}$ is a stationary Markov process on $\{-1, 1\}^{\mathbb{Z}^d}$ with marginal distribution $\mu^{+,\beta,h}$ governed by flip rate intensities

$$(10) \qquad C(s, \sigma) = \exp\left(-\beta \sum_{t \in \mathbb{Z}^d : \{t,s\} \in E} \sigma(t)\sigma(s) - h\sigma(s)\right);$$

see [19], Theorem 2.7, page 139. It is also possible to construct $\Psi^{+,\beta,h}$ directly on $\{-1, 1\}^{\mathbb{Z}^d}$ without going through the limiting procedure. Furthermore, there are several possible choices of flip rate intensities that can be used to construct a stationary and reversible Markov process on $\{-1, 1\}^{\mathbb{Z}^d}$ with marginal distribution $\mu^{+,\beta,h}$. In [19], a stochastic Ising model is defined to be any spin system with flip rate intensities $\{C(s, \sigma)\}_{s \in \mathbb{Z}^d, \sigma \in \{-1,1\}^{\mathbb{Z}^d}}$ satisfying that, for each $s \in \mathbb{Z}^d$,

$$(11) \qquad C(s, \sigma) \exp\left(\beta \sum_{\substack{\{t,s\} \in E \\ t \in \mathbb{Z}^d}} \sigma(t)\sigma(s) + h\sigma(s)\right)$$

is independent of $\sigma(s)$. Therefore, when we refer to a stochastic Ising model $\Psi^{+,\beta,h}$ with marginal distribution $\mu^{+,\beta,h}$, we will have this definition in mind. It is particularly easy to see that (11) (or the condition of detailed balance as it is often referred to) is satisfied for the flip rate intensities of (10), but there are many other rates satisfying this. It is known that the set of so-called Gibbs states are exactly the same as the class of reversible measures with respect to the flip rates satisfying (11); see [19], pages 190–196. Note also that the condition specified in (11) with $\mathbb{Z}^d$ replaced by $\Lambda_n$ is equivalent to that of (9) (modified with the boundary condition removed).

While we defined above stochastic Ising models on $\{-1,1\}^{\mathbb{Z}^d}$, this construction can be done on more general graphs (see [19]).

2.4. *The random cluster model.* Unlike all other models in this paper, the random cluster model deals with configurations on the edges $E$ of a graph $G = (S, E)$. We will review the definition of the regular random cluster measure on general finite graphs and the "wired" random cluster measure on $\Lambda_n \subseteq \mathbb{Z}^d$. We will also recall the limiting measures and in the next subsection the connection between the random cluster model and the Ising model. In doing so we will follow the outlines of [9] and [13] closely.

Take a finite graph $G = (S, E)$. Define the random cluster measure $\nu_G^{p,q}$ on $\{0,1\}^E$ with parameters $p \in [0,1]$ and $q > 0$ as the probability measure which assigns to the configuration $\eta \in \{0,1\}^E$ the probability

$$(12) \qquad \nu_G^{p,q}(\eta) = \frac{q^{k(\eta)}}{Z} \prod_{e \in E} p^{\eta(e)} (1-p)^{1-\eta(e)}.$$

Here $Z$ is again a normalization constant and $k(\eta)$ is the number of connected components of $\eta$. From now on we will always take $q = 2$ and, therefore, we will suppress q in the notation.

Take $G_n = (S_n, E_n)$, where $S_n = \Lambda_{n+1} \subseteq \mathbb{Z}^d$ and $E_n$ is the set of all nearest neighbor pairs of $\Lambda_{n+1}$. Write $\nu_n^p$ for $\nu_{G_n}^p$, and define

$$(13) \qquad \tilde{\nu}_n^p(\cdot) = \nu_n^p(\cdot | \text{all edges of } E_n \text{ with both}$$

$$\text{end sites in } \Lambda_{n+1} \setminus \Lambda_n \text{ are present}).$$

This is the so-called "wired" random cluster measure. It is called "wired" since all edges of the boundary are present. It is immediate from the defining equations (12) and (13), that, for $e \in E_n$ and any $\xi \in \{0,1\}^{E_n \setminus e}$,

$$
\begin{aligned}
(14) \quad & \tilde{\nu}_n^p(\eta(e) = 1 | \eta(E_n \setminus e) = \xi) \\
& = \begin{cases} p, & \text{if the endpoints of } e \text{ are connected in } \xi, \\ \dfrac{p}{2-p}, & \text{otherwise.} \end{cases}
\end{aligned}
$$

One can show (see [9] or [13]) that when n tends to infinity, the probability measures $\{\tilde{\nu}_n^p\}_{n \in \mathbb{N}^+}$ converge to a probability measure $\tilde{\nu}^p$. Furthermore, the construction of $\tilde{\nu}_n^p$ on $\{0,1\}^{E_n}$ can be done on any finite subgraph by connecting all sites of the boundary of the graph with each other. As a consequence, we can also define random cluster measures on more general graphs than $\mathbb{Z}^d$; see, for example, [11].

2.5. *The random cluster model and the Ising model.* Take $G_n = (S_n, E_n)$ as in Section 2.4. As in [13], let $\mathbf{P}_n^p$ be the probability measure on $\{-1,1\}^{S_n} \times \{0,1\}^{E_n}$ defined in the following way:

1. Assign each site of $\Lambda_{n+1} \setminus \Lambda_n$ and every edge with both endpoints in $\Lambda_{n+1} \setminus \Lambda_n$ the value 1.
2. Assign each site of $\Lambda_n$ the value 1 or $-1$ with equal probability, assign each edge with not more than one endpoint in $\Lambda_{n+1} \setminus \Lambda_n$ the value 0 or 1 with probabilities $1-p$ and $p$, respectively. Do this independently for all sites and edges.
3. Condition on the event that no two sites with different spins have an open edge connecting them.

One can then check that $\mathbf{P}_n^p(\sigma, \{0,1\}^{E_n}) = \mu_n^{+,\beta}(\sigma)$ with $\beta = -\log(1-p)/2$, and that $\mathbf{P}_n^p(\{-1,1\}^{S_n}, \eta) = \tilde{\nu}_n^p(\eta)$. Here, $\mathbf{P}_n^p(\sigma, \{0,1\}^{E_n})$ is just the marginal in the first coordinate of $\mathbf{P}_n^p$. The same kind of construction can be carried out on any finite graph $G = (S, E)$.

2.6. *The contact process.* Consider a graph $G = (S, E)$ of bounded degree. In the contact process the state space is $\{0,1\}^S$. Let $\lambda > 0$, and define the flip rate intensities to be

$$
C(s, \sigma) = \begin{cases} 1, & \text{if } \sigma(s) = 1, \\ \lambda \displaystyle\sum_{(s',s) \in E} \sigma(s'), & \text{if } \sigma(s) = 0. \end{cases}
$$

If we let the initial distribution be $\sigma \equiv 1$, the distribution of this process at time $t$, which we will denote by $\delta_1 T_\lambda(t)$, is known to converge as $t$ tends to infinity. This is simply because it is a so-called "attractive" process and $\sigma \equiv 1$ is the maximal state and $\{\delta_1 T_\lambda(t)\}$ is stochastically decreasing; see [19], page 265. This limiting distribution will be referred to as the upper invariant measure for the contact process with parameter $\lambda$ and will be denoted by $\mu_\lambda$. We then let $\Psi^\lambda$ denote the stationary Markov process on $\{0,1\}^S$ with initial (and invariant) distribution $\mu_\lambda$.

**3. $\varepsilon$-movability for monotone measures.** In this section we prove movability results for classes of monotone measures. The finite case is covered by Lemma 3.3, while the countable case is discussed in Proposition 3.4. In this section we will always assume that our measures have full support.

For any $|S| < \infty$, $s \in S$, $\xi \in \{0,1\}^{S \setminus s}$ and probability measure $\mu$ on $\{0,1\}^S$, write $\mu^{(*,\varepsilon)}(i|\xi)$ for $\mu^{(*,\varepsilon)}(\sigma(s) = i | \sigma(S \setminus s) = \xi)$, $\mu^{(*,\varepsilon)}(i \cap \xi)$ for $\mu^{(*,\varepsilon)}(\{\sigma(s) = i\} \cap \{\sigma(S \setminus s) = \xi\})$ and $\mu^{(*,\varepsilon)}(\xi)$ for $\mu^{(*,\varepsilon)}(\sigma(S \setminus s) = \xi)$. Here "$*$" can represent either $+$ or $-$ and $i \in \{0,1\}$. Note that $s$ is suppressed in the notation and so should be understood from context.

We begin with an easy lemma whose proof is left to the reader. The idea is that if the configuration outside of $s$ is $\xi$ under $\mu^{(-,\varepsilon)}$, it must have been at least as large under $\mu$ "before flipping some 1's to 0's"; then use monotonicity.

LEMMA 3.1. *Assume that $\mu$ is a monotone probability measure on $\{0,1\}^S$, where $|S| < \infty$. Take $s \in S$ and let $\xi \in \{0,1\}^{S \setminus s}$. Then, for any $\varepsilon > 0$, we have that*

$$\mu^{(-,\varepsilon)}(1|\xi) \geq (1-\varepsilon)\mu(1|\xi)$$

*and that*

$$\mu^{(+,\varepsilon)}(0|\xi) \geq (1-\varepsilon)\mu(0|\xi).$$

The next lemma will be used to prove Lemma 3.3.

LEMMA 3.2. *Assume that $\mu$ is a monotone probability measure on $\{0,1\}^S$, where $|S| < \infty$. For any $\varepsilon > 0$, $\mu^{(-,\varepsilon)}$ is also monotone.*

PROOF. Let $s \in S$ be arbitrary, $X \sim \mu$ and $X^{(-,\varepsilon)} \sim \mu^{(-,\varepsilon)}$. For any $\delta, \eta \in \{0,1\}^{S \setminus s}$, define the probability measures $\mu_\delta$ and $\mu_\eta$ on $\{0,1\}^{S \setminus s}$ by letting $\mu_\delta(\mathcal{A}) = \mathbb{P}(X \in \mathcal{A} | X^{(-,\varepsilon)}(S \setminus s) \equiv \delta)$ and $\mu_\eta(\mathcal{A}) = \mathbb{P}(X \in \mathcal{A} | X^{(-,\varepsilon)}(S \setminus s) \equiv \eta)$ for every event $\mathcal{A}$ in $\{0,1\}^{S \setminus s}$, respectively. We will prove that

$$(15) \qquad\qquad \mu_\delta \preceq \mu_\eta \qquad \forall \delta \preceq \eta.$$

This will give us [since $\mathbb{P}(X(s) = 1 | X(S \setminus s) \equiv \eta)$ is an increasing function of $\eta$] that

$$\mathbb{P}(X^{(-,\varepsilon)}(s) = 1 | X^{(-,\varepsilon)}(S \setminus s) \equiv \eta)$$

$$= (1-\varepsilon) \int_{\tilde{\eta} \in \{0,1\}^{S \setminus s}} \mathbb{P}(X(s) = 1 | X(S \setminus s) \equiv \tilde{\eta}) \, d\mu_\eta(\tilde{\eta})$$

$$\geq (1-\varepsilon) \int_{\tilde{\eta} \in \{0,1\}^{S \setminus s}} \mathbb{P}(X(s) = 1 | X(S \setminus s) \equiv \tilde{\eta}) \, d\mu_\delta(\tilde{\eta})$$

$$= \mathbb{P}(X^{(-,\varepsilon)}(s) = 1 | X^{(-,\varepsilon)}(S \setminus s) \equiv \delta).$$

Since $s$ was chosen arbitrarily, this would prove the statement.

We now prove (15). Define for $\eta \preceq \tilde{\eta}$ $d(\tilde{\eta}, \eta) := |\{t \in S \setminus s : \tilde{\eta}(t) = 1\}| - |\{t \in S \setminus s : \eta(t) = 1\}|$ and $d(\tilde{\eta}, 0) = |\{t \in S \setminus s : \tilde{\eta}(t) = 1\}|$. Here $|\cdot|$ denotes cardinality. Let $\mu_{S \setminus s}(\eta) = \mathbb{P}(X(S \setminus s) \equiv \eta)$ and define $\mu_{S \setminus s}^{(-,\varepsilon)}$ similarly. We have that, for $\eta \preceq \tilde{\eta}$,

$$\mu_\eta(\tilde{\eta}) = \mathbb{P}(X^{(-,\varepsilon)}(S \setminus s) \equiv \eta | X(S \setminus s) \equiv \tilde{\eta}) \frac{\mu_{S \setminus s}(\tilde{\eta})}{\mu_{S \setminus s}^{(-,\varepsilon)}(\eta)} \tag{16}$$

$$= \varepsilon^{d(\tilde{\eta}, \eta)}(1-\varepsilon)^{d(\eta, 0)} \frac{\mu_{S \setminus s}(\tilde{\eta})}{\mu_{S \setminus s}^{(-,\varepsilon)}(\eta)}. \tag{17}$$

It is well known that $\mu$ being monotone implies that, for every, $\tilde{\delta}$, $\tilde{\eta}$,

$$\mu_{S \setminus s}(\tilde{\eta} \vee \tilde{\delta}) \mu_{S \setminus s}(\tilde{\eta} \wedge \tilde{\delta}) \geq \mu_{S \setminus s}(\tilde{\eta}) \mu_{S \setminus s}(\tilde{\delta}). \tag{18}$$

By a simple modification of Theorem 2.9, page 75 of [19], it is enough for us to show that

$$\mu_\eta(\tilde{\eta} \vee \tilde{\delta}) \mu_\delta(\tilde{\eta} \wedge \tilde{\delta}) \geq \mu_\eta(\tilde{\eta}) \mu_\delta(\tilde{\delta}) \tag{19}$$

for all $\tilde{\eta} \succeq \eta$ and $\tilde{\delta} \succeq \delta$ to show (15). An elementary calculation shows that

$$d(\tilde{\eta} \vee \tilde{\delta}, \eta) + d(\tilde{\eta} \wedge \tilde{\delta}, \delta) = d(\tilde{\eta}, \eta) + d(\tilde{\delta}, \delta). \tag{20}$$

We therefore get

$$\mu_\eta(\tilde{\eta} \vee \tilde{\delta}) \mu_\delta(\tilde{\eta} \wedge \tilde{\delta})$$

$$= \varepsilon^{d(\tilde{\eta} \vee \tilde{\delta}, \eta) + d(\tilde{\eta} \wedge \tilde{\delta}, \delta)} (1-\varepsilon)^{d(\eta, 0) + d(\delta, 0)} \frac{\mu_{S \setminus s}(\tilde{\eta} \vee \tilde{\delta})}{\mu_{S \setminus s}^{(-,\varepsilon)}(\eta)} \frac{\mu_{S \setminus s}(\tilde{\eta} \wedge \tilde{\delta})}{\mu_{S \setminus s}^{(-,\varepsilon)}(\delta)}$$

$$\geq \varepsilon^{d(\tilde{\eta}, \eta) + d(\tilde{\delta}, \delta)} (1-\varepsilon)^{d(\eta, 0) + d(\delta, 0)} \frac{\mu_{S \setminus s}(\tilde{\eta})}{\mu_{S \setminus s}^{(-,\varepsilon)}(\eta)} \frac{\mu_{S \setminus s}(\tilde{\delta})}{\mu_{S \setminus s}^{(-,\varepsilon)}(\delta)} = \mu_\eta(\tilde{\eta}) \mu_\delta(\tilde{\delta}),$$

where (16) is used in the first and last equality and equations (18) and (20) are used in the inequality.   □

LEMMA 3.3.  *Let $\mu_1, \mu_2$ be probability measures on $\{0,1\}^S$, where $|S| < \infty$. Assume that $\mu_2$ is monotone and that*

$$A := \inf_{\substack{s \in S \\ \xi \in \{0,1\}^{S \setminus s}}} [\mu_2(\sigma(s) = 1 | \sigma(S \setminus s) \equiv \xi) - \mu_1(\sigma(s) = 1 | \sigma(S \setminus s) \equiv \xi)] > 0.$$

*Then for any choice of $\varepsilon > 0$, such that*

$$A > \frac{1}{1-\varepsilon} - 1,$$

*we have*

$$\mu_1 \preceq \mu_2^{(-,\varepsilon)}.$$

*Hence, $(\mu_1, \mu_2)$ is downward movable.*

PROOF. Monotonicity of $\mu_2$, Lemma 3.1, the definition of $A$ and our choice of $\varepsilon$ give us that, for any $s \in S$ and $\xi \in \{0,1\}^{S \setminus s}$,

$$\mu_2^{(-,\varepsilon)}(1|\xi) \geq (1-\varepsilon)\mu_2(1|\xi) \geq (1-\varepsilon)(A + \mu_1(1|\xi))$$

$$\geq (1-\varepsilon)\frac{\mu_1(1|\xi)}{1-\varepsilon} = \mu_1(1|\xi).$$

By Lemma 3.2, $\mu_2^{(-,\varepsilon)}$ is monotone and so $\forall \tilde{\xi} \preceq \xi$,

$$\mu_1(1|\tilde{\xi}) \leq \mu_2^{(-,\varepsilon)}(1|\tilde{\xi}) \leq \mu_2^{(-,\varepsilon)}(1|\xi).$$

The proof is completed by the use of Holley's inequality, Theorem 2.1. □

PROPOSITION 3.4. *Let $S$ be any finite or countable set and consider $(S_n)_{n \in \mathbb{N}^+}$, a collection of sets such that $|S_n| < \infty$ $\forall n \in \mathbb{N}^+$ and $S_n \uparrow S$. Let $(\mu_{1,n})_{n \in \mathbb{N}^+}$, $(\mu_{2,n})_{n \in \mathbb{N}^+}$, be two collections of probability measures, where $\mu_{1,n}, \mu_{2,n}$ are probability measures on $\{0,1\}^{S_n}$ for every $n \in \mathbb{N}^+$. Furthermore, assume that all of the probability measures $(\mu_{1,n})_{n \in \mathbb{N}^+}((\mu_{2,n})_{n \in \mathbb{N}^+})$ are monotone, that $\mu_{1,n} \to \mu_1$ and that $\mu_{2,n} \to \mu_2$. Set*

$$A_n := \inf_{\substack{s \in S_n \\ \xi \in \{0,1\}^{S_n \setminus s}}} [\mu_{2,n}(\sigma(s) = 1|\sigma(S \setminus s) \equiv \xi) - \mu_{1,n}(\sigma(s) = 1|\sigma(S \setminus s) \equiv \xi)].$$

*If*

$$\inf_{n \in \mathbb{N}^+} A_n > 0,$$

*then $(\mu_1, \mu_2)$ is both upward and downward movable.*

PROOF. Take $\varepsilon > 0$ such that

$$\inf_{n \in \mathbb{N}^+} A_n > \frac{1}{1-\varepsilon} - 1.$$

With this choice of $\varepsilon$, Lemma 3.3 says that $(\mu_{1,n}, \mu_{2,n})$ is upward (downward) $\varepsilon$-movable. Since $\mu_{1,n} \to \mu_1$ and $\mu_{2,n} \to \mu_2$, we easily get that $\mu_{2,n}^{(-,\varepsilon)} \to \mu_2^{(-,\varepsilon)}$ and $\mu_{1,n}^{(+,\varepsilon)} \to \mu_1^{(+,\varepsilon)}$. Furthermore, since the relations

$$\mu_{1,n} \preceq \mu_{2,n}^{(-,\varepsilon)}$$

and

$$\mu_{1,n}^{(+,\varepsilon)} \preceq \mu_{2,n}$$

are easily seen to be preserved under weak limits, we get that

$$\mu_1 \preceq \mu_2^{(-,\varepsilon)} \quad \text{and} \quad \mu_1^{(+,\varepsilon)} \preceq \mu_2. \qquad\qquad □$$

**4. $\varepsilon$-movability for the contact process and a 0–1 law.** The conditions in our next proposition might seem overly technical; however, these represent the essential features of the contact process (after a small suitable time rescaling) and, therefore, we feel it is instructive to highlight these features. In Proposition 4.1 and Lemmas 5.1, 5.2 and 8.1 we will use the so-called graphical representation to define our processes; see, for instance, [19], page 172.

PROPOSITION 4.1. *Let $\mu_1$ and $\mu_2$ be two probability measures defined on $\{0,1\}^S$, where $S$ is a countable set. Assume that $\mu_1 \preceq \mu_2$ and that there exists two stationary Markov processes $\Psi_1$ and $\Psi_2$, governed by flip rate intensities $\{C_1(s,\sigma_1)\}_{s \in S, \sigma_1 \in \{0,1\}^S}$ and $\{C_2(s,\sigma_2)\}_{s \in S, \sigma_2 \in \{0,1\}^S}$, respectively, and with marginal distributions $\mu_1$ and $\mu_2$. Assume that $C_1 \preceq C_2$ [conditions (7) and (8) of the Introduction]. Consider the following conditions:*

1. *There exists an $\varepsilon_1 > 0$ such that*

$$C_2(s,\sigma_2) - C_1(s,\sigma_1) \geq \varepsilon_1$$

(21)

$$\forall s \in S, \ \forall \sigma_2 \succeq \sigma_1 \ s.t. \ \sigma_2(s) = 0 \ and \ C_1(s,\sigma_1) \neq 0.$$

2. *There exists an $\varepsilon_2 > 0$ such that*

$$C_1(s,\sigma_1) - C_2(s,\sigma_2) \geq \varepsilon_2$$

(22)

$$\forall s \in S, \ \forall \sigma_2 \succeq \sigma_1 \ s.t. \ \sigma_1(s) = 1 \ and \ C_2(s,\sigma_2) \neq 0.$$

3. *There exists an $\varepsilon_3 > 0$ such that*

(23)          $$C_1(s,\sigma_1) \geq \varepsilon_3 \qquad \forall s \in S, \ \forall \sigma_1 \ s.t. \ \sigma_1(s) = 1.$$

4. *There exists an $\varepsilon_4 > 0$ such that*

(24)          $$C_2(s,\sigma_2) \geq \varepsilon_4 \qquad \forall s \in S, \ \forall \sigma_2 \ s.t. \ \sigma_2(s) = 0.$$

*If conditions 1, 2 and 3 are satisfied, then $(\mu_1, \mu_2)$ is downward movable. If conditions 1, 2 and 4 are satisfied, then $(\mu_1, \mu_2)$ is upward movable.*

PROOF. We will prove the first statement, the second follows by symmetry. Define

$$\lambda := \sup_{s,\sigma_2 \, : \, \sigma_2(s)=0} C_2(s,\sigma_2) + \sup_{s,\sigma_1 \, : \, \sigma_1(s)=1} C_1(s,\sigma_1).$$

Our aim is to construct a coupling of the processes $\{X_{1,t}\}_{t \geq 0} \sim \Psi_1$ and $\{X_{2,t}\}_{t \geq 0} \sim \Psi_2$ such that $X_{1,t} \preceq X_{2,t} \ \forall t \geq 0$ in such a way that we prove the proposition. Before presenting the actual coupling, we will discuss the idea behind it. For every site $s \in S$, associate an independent Poisson process with parameter $\lambda$. Next, let $\{U_{s,k}\}_{s \in S, k \geq 1}$ and $\{U'_{s,k}\}_{s \in S, k \geq 1}$ be independent uniform $[0, 1]$ random variables also independent of the Poisson processes. If $\tau$ is an arrival time for the Poisson process at site $s$, we write $U_{s,\tau}$ for $U_{s,k}$, where $k$ is such that $\tau$ is the $k$th arrival of the Poisson process at site $s$. Now, let $\tau$ be an arrival time for the Poisson process associated to a site $s$. For $i \in \{1, 2\}$, let $X_{i,\tau^-}$ and $X_{i,\tau^+}$ denote the configurations before and after the arrival. We will let the outcome of $U_{s,\tau}$ decide what happens with the $\{X_{2,t}\}_{t \geq 0}$ process at time $t = \tau$, and then we will let $U'_{s,\tau}$, together with $U_{s,\tau}$, decide what happens with the $\{X_{1,t}\}_{t \geq 0}$ process at time $t = \tau$. As we will see, we will do this so that $X_{1,t} \preceq X_{2,t}$ for all $t \geq 0$. Furthermore, we will do this in such a way that there exists an $\varepsilon \in (0, 1)$ such that if $U'_{s,\tau} \geq 1 - \varepsilon$, then $X_{1,\tau^+}(s) = 0$ regardless of the outcome of $U_{s,\tau}$. Consider now the process $\{X_t^\varepsilon\}_{t \geq 0}$ we get by taking $X_0^\varepsilon(s) = 1$ for every $s \in S$ and letting $\{X_t^\varepsilon(s)\}_{t \geq 0}$ be updated at every arrival time $\tau$ for the Poisson process associated to $s$, and updated in such a way that $X_{\tau^+}^\varepsilon(s) = 0$ if $U'_{s,\tau} \geq 1 - \varepsilon$, and $X_{\tau^+}^\varepsilon(s) = 1$ if $U'_{s,\tau} < 1 - \varepsilon$. Of course, the distribution of $X_t^\varepsilon$ will converge to $\pi_{1-\varepsilon}$. Observe that whenever $X_t^\varepsilon(s) = 0$, we have that $X_{1,t}(s) = 0$. Therefore, we can conclude that

$$(25) \qquad X_{1,t} \preceq \min(X_{2,t}, X_t^\varepsilon) \qquad \forall t \geq 0.$$

Furthermore, since the process $\{X_t^\varepsilon\}_{t \geq 0}$ does not depend on any $U_{s,\tau}$, we have that $X_t^\varepsilon(s)$ is conditionally independent of $X_{2,t}$ if there has been an arrival for the Poisson process associated to $s$ before time $t$. Let $s_i$, $i \in \{1, \ldots, n\}$, be distinct sites in $S$ and let $\mathcal{A}_t$ be the event that all Poisson processes associated to $s_1$ through $s_n$ have had an arrival by time $t$. Of course, $\mathbb{P}(\mathcal{A}_t) = (1 - e^{-\lambda t})^n$ and so we get that

$$\begin{aligned}
&\mathbb{P}(X_{2,t} X_t^\varepsilon(s_1) = \cdots = X_{2,t} X_t^\varepsilon(s_n) = 1) \\
&\quad = \mathbb{P}(X_{2,t} X_t^\varepsilon(s_1) = \cdots = X_{2,t} X_t^\varepsilon(s_n) = 1 | \mathcal{A}_t) \mathbb{P}(\mathcal{A}_t) \\
&\qquad + \mathbb{P}(X_{2,t} X_t^\varepsilon(s_1) = \cdots = X_{2,t} X_t^\varepsilon(s_n) = 1 | \mathcal{A}_t^c) \mathbb{P}(\mathcal{A}_t^c) \\
&\quad = \mathbb{P}(X_{2,t}(s_1) = \cdots = X_{2,t}(s_n) = 1 | \mathcal{A}_t) \\
&\qquad \times \mathbb{P}(X_t^\varepsilon(s_1) = \cdots = X_t^\varepsilon(s_n) = 1 | \mathcal{A}_t) \mathbb{P}(\mathcal{A}_t) \\
&\qquad + \mathbb{P}(X_{2,t} X_t^\varepsilon(s_1) = \cdots = X_{2,t} X_t^\varepsilon(s_n) = 1 | \mathcal{A}_t^c) \mathbb{P}(\mathcal{A}_t^c) \\
&\quad = \mathbb{P}(X_{2,t}(s_1) = \cdots = X_{2,t}(s_n) = 1 | \mathcal{A}_t) \mathbb{P}(\mathcal{A}_t)(1 - \varepsilon)^n \\
&\qquad + \mathbb{P}(X_{2,t} X_t^\varepsilon(s_1) = \cdots = X_{2,t} X_t^\varepsilon(s_n) = 1 | \mathcal{A}_t^c) \mathbb{P}(\mathcal{A}_t^c) \\
&\quad = \mathbb{P}(\{X_{2,t}(s_1) = \cdots = X_{2,t}(s_n) = 1\} \cap \mathcal{A}_t)(1 - \varepsilon)^n \\
&\qquad + \mathbb{P}(X_{2,t} X_t^\varepsilon(s_1) = \cdots = X_{2,t} X_t^\varepsilon(s_n) = 1 | \mathcal{A}_t^c) \mathbb{P}(\mathcal{A}_t^c) \\
&\quad \geq (\mathbb{P}(X_{2,t}(s_1) = \cdots = X_{2,t}(s_n) = 1) - \mathbb{P}(\mathcal{A}_t^c))(1 - \varepsilon)^n \\
&\qquad + \mathbb{P}(X_{2,t} X_t^\varepsilon(s_1) = \cdots = X_{2,t} X_t^\varepsilon(s_n) = 1 | \mathcal{A}_t^c) \mathbb{P}(\mathcal{A}_t^c)
\end{aligned}$$

$$= \mathbb{P}(X_{2,t}(s_1) = \cdots = X_{2,t}(s_n) = 1)(1-\varepsilon)^n$$
$$+ \mathbb{P}(\mathcal{A}_t^c)(\mathbb{P}(X_{2,t}X_t^\varepsilon(s_1) = \cdots = X_{2,t}X_t^\varepsilon(s_n) = 1|\mathcal{A}_t^c) - (1-\varepsilon)^n)$$
$$= \mu_2^{(-,\varepsilon)}(\sigma(s_1) = \cdots = \sigma(s_n) = 1)$$
$$+ \mathbb{P}(\mathcal{A}_t^c)(\mathbb{P}(X_{2,t}X_t^\varepsilon(s_1) = \cdots = X_{2,t}X_t^\varepsilon(s_n) = 1|\mathcal{A}_t^c) - (1-\varepsilon)^n)$$
$$\overset{t \to \infty}{\longrightarrow} \mu_2^{(-,\varepsilon)}(\sigma(s_1) = \cdots = \sigma(s_n) = 1).$$

In addition,

$$\mathbb{P}(X_{2,t}(s_1) = \cdots = X_{2,t}(s_n) = 1 \cap \mathcal{A}_t)(1-\varepsilon)^n$$
$$\leq \mathbb{P}(X_{2,t}(s_1) = \cdots = X_{2,t}(s_n) = 1)(1-\varepsilon)^n$$
$$= \mu_2^{(-,\varepsilon)}(\sigma(s_1) = \cdots = \sigma(s_n) = 1).$$

Hence, by inclusion exclusion, we have that the distribution of $\min(X_{2,t}, X_t^\varepsilon)$ approaches $\mu_2^{(-,\varepsilon)}$ as $t$ tends to infinity. So by first taking the limit in (25), we get that $\mu_1 \preceq \mu_2^{(-,\varepsilon)}$, as desired.

Now to the construction. Take $X_{1,0} \sim \mu_1$, $X_{2,0} \sim \mu_2$, such that $X_{1,0} \preceq X_{2,0}$. Let $\tau$ be an arrival time for the Poisson process associated to $s$. Take $U_{s,\tau}$ and $U'_{s,\tau}$. The following transition rules apply:

| $X_{2,\tau^-}$ | $X_{2,\tau^+}$ | if |
|---|---|---|
| 0 | 1 | $U_{s,\tau} \leq \dfrac{C_2(s, X_{2,\tau^-})}{\lambda}$ |
| 1 | 0 | $U_{s,\tau} \geq \dfrac{\lambda - C_2(s, X_{2,\tau^-})}{\lambda}.$ |

It is easy to check that the process $\{X_{2,t}\}_{t \geq 0}$ thus constructed will have the right flip-rate intensities. The construction of $\{X_{1,t}\}_{t \geq 0}$ is slightly more complicated. If $C_2(s, X_{2,\tau^-}) = 0$ and $X_{2,\tau^-}(s) = 0$, then it follows from (7) that $C_1(s, X_{1,\tau^-}) = 0$, and in that case we interpret $\frac{C_1(s, X_{1,\tau^-})}{C_2(s, X_{2,\tau^-})}$ as 0. Observe that $C_2(s, X_{2,\tau^-})$ can be 0 when $X_{2,\tau^-}(s) = 1$, but it will not cause any problems. With these observations in mind, these are the transition rules we

apply:

| $(X_{1,\tau^-}, X_{2,\tau^-})$ | $(X_{1,\tau^+}, X_{2,\tau^+})$ | if |
|---|---|---|
| $(0,0)$ | $(1,1)$ | $U_{s,\tau} \leq \dfrac{C_2(s, X_{2,\tau^-})}{\lambda}$ and $U'_{s,\tau} \leq \dfrac{C_1(s, X_{1,\tau^-})}{C_2(s, X_{2,\tau^-})}$ |
| $(0,0)$ | $(0,1)$ | $U_{s,\tau} \leq \dfrac{C_2(s, X_{2,\tau^-})}{\lambda}$ and $U'_{s,\tau} > \dfrac{C_1(s, X_{1,\tau^-})}{C_2(s, X_{2,\tau^-})}$ |
| $(0,0)$ | $(0,0)$ | otherwise |
| $(0,1)$ | $(0,0)$ | $U_{s,\tau} \geq \dfrac{\lambda - C_2(s, X_{2,\tau^-})}{\lambda}$ |
| $(0,1)$ | $(1,1)$ | $U_{s,\tau} < \dfrac{\sup_{s,\sigma_2 : \sigma_2(s)=0} C_2(s,\sigma_2)}{\lambda}$ and $U'_{s,\tau} \leq \dfrac{C_1(s, X_{1,\tau^-})}{\sup_{s,\sigma_2 : \sigma_2(s)=0} C_2(s,\sigma_2)}$ |
| $(0,1)$ | $(0,1)$ | otherwise |
| $(1,1)$ | $(0,0)$ | $U_{s,\tau} \geq \dfrac{\lambda - C_2(s, X_{2,\tau^-})}{\lambda}$ |
| $(1,1)$ | $(0,1)$ | $U_{s,\tau} < \dfrac{\lambda - C_2(s, X_{2,\tau^-})}{\lambda}$ and $U'_{s,\tau} \geq \dfrac{\lambda - C_1(s, X_{1,\tau^-})}{\lambda - C_2(s, X_{2,\tau^-})}$ |
| $(1,1)$ | $(1,1)$ | otherwise. |

It is not difficult to check that all flip rate intensities are correct and that $X_{1,t} \preceq X_{2,t}$ for all $t \geq 0$. Observe that, by the definition of $\lambda$, the events $\{U_{s,\tau} \geq \frac{\lambda - C_2(s, X_{2,\tau^-})}{\lambda}\}$ and $\{U_{s,\tau} < \frac{\sup_{s,\sigma_2 : \sigma_2(s)=0} C_2(s,\sigma_2)}{\lambda}\}$ are disjoint when $(X_{1,\tau^-}, X_{2,\tau^-}) = (0,1)$.

We now want to show that there exists an $\varepsilon > 0$ so that $U'_{s,\tau} \geq 1 - \varepsilon$ implies that $X_{1,\tau^+}(s) = 0$. Note that if $(X_{1,\tau^-}, X_{2,\tau^-}) = (0,0)$ and $C_1(s, X_{1,\tau^-}) > 0$ [$\Rightarrow C_2(s, X_{2,\tau^-}) > 0$], then

$$\frac{C_1(s, X_{1,\tau^-})}{C_2(s, X_{2,\tau^-})} \leq \frac{C_2(s, X_{2,\tau^-}) - \varepsilon_1}{C_2(s, X_{2,\tau^-})} \leq 1 - \frac{\varepsilon_1}{\sup_{s,\sigma_2 : \sigma_2(s)=0} C_2(s,\sigma_2)} < 1$$

and if $(X_{1,\tau^-}, X_{2,\tau^-}) = (0,0)$ and $C_1(s, X_{1,\tau^-}) = 0$, then

$$\frac{C_1(s, X_{1,\tau^-})}{C_2(s, X_{2,\tau^-})} = 0.$$

Furthermore, if $(X_{1,\tau^-}, X_{2,\tau^-}) = (0,1)$ and $C_1(s, X_{1,\tau^-}) > 0$, then

$$\frac{C_1(s, X_{1,\tau^-})}{\sup_{s,\sigma_2 : \sigma_2(s)=0} C_2(s,\sigma_2)} \leq 1 - \frac{\varepsilon_1}{\sup_{s,\sigma_2 : \sigma_2(s)=0} C_2(s,\sigma_2)} < 1,$$

while again if $(X_{1,\tau^-}, X_{2,\tau^-}) = (0,1)$ and $C_1(s, X_{1,\tau^-}) = 0$, then the 0 never changes to a 1. Finally, if $(X_{1,\tau^-}, X_{2,\tau^-}) = (1,1)$ and $C_2(s, X_{2,\tau^-}) > 0$ [$\Rightarrow C_1(s, X_{1,\tau^-}) > 0$], then

$$\frac{\lambda - C_1(s, X_{1,\tau^-})}{\lambda - C_2(s, X_{2,\tau^-})} \leq \frac{\lambda - C_2(s, X_{2,\tau^-}) - \varepsilon_2}{\lambda - C_2(s, X_{2,\tau^-})}$$

$$\leq 1 - \frac{\varepsilon_2}{\lambda - C_2(s, X_{2,\tau^-})}$$

$$\leq 1 - \frac{\varepsilon_2}{\lambda},$$

and if $(X_{1,\tau^-}, X_{2,\tau^-}) = (1,1)$ and $C_2(s, X_{2,\tau^-}) = 0$,

$$\frac{\lambda - C_1(s, X_{1,\tau^-})}{\lambda - C_2(s, X_{2,\tau^-})} \leq \frac{\lambda - \varepsilon_3}{\lambda} = 1 - \frac{\varepsilon_3}{\lambda} < 1.$$

Therefore, whenever

$$U'_{s,\tau} \geq \max\left(1 - \frac{\varepsilon_1}{\sup_{s, \sigma_2 : \sigma_2(s) = 0} C_2(s, \sigma_2)}, 1 - \frac{\varepsilon_2}{\lambda}, 1 - \frac{\varepsilon_3}{\lambda}\right),$$

we have that $X_{1,\tau^+}(s) = 0$ regardless of the outcome of $U_{s,\tau}$. Therefore, $(\mu_1, \mu_2)$ is downward $\varepsilon$-movable where

$$\varepsilon := 1 - \max\left(1 - \frac{\varepsilon_1}{\sup_{s, \sigma_2 : \sigma_2(s) = 0} C_2(s, \sigma_2)}, 1 - \frac{\varepsilon_2}{\lambda}, 1 - \frac{\varepsilon_3}{\lambda}\right)$$

$$= \min\left(\frac{\varepsilon_1}{\sup_{s, \sigma_2 : \sigma_2(s) = 0} C_2(s, \sigma_2)}, \frac{\varepsilon_2}{\lambda}, \frac{\varepsilon_3}{\lambda}\right). \qquad \square$$

PROOF OF THEOREM 1.11.    Take $\delta > 0$ such that $\lambda_1(1 + \delta) < \lambda_2$ and consider the process $\{X_t\}_{t \geq 0}$ constructed in the following way. Take $X_0 \equiv 1$ and let the process evolve with flip rate intensities

(26)    $$C_1(s, \sigma) = \begin{cases} 1 + \delta, & \text{if } \sigma(s) = 1, \\ \lambda_1(1 + \delta) \sum_{s' \sim s} \sigma(s'), & \text{if } \sigma(s) = 0. \end{cases}$$

Denote the limiting distribution of $X_t$ as $t$ tends to infinity by $\mu_{1+\delta, \lambda_1(1+\delta)}$. It is easy to see that this process is just a time-scaling of the contact process constructed in Section 2.6 with parameter $\lambda_1$. Recall that that process had limiting distribution $\mu_{\lambda_1}$, the upper invariant measure for the contact process. Thus, we have $\mu_{\lambda_1} = \mu_{1+\delta, \lambda_1(1+\delta)}$. By Proposition 4.1 with $C_1$ as above and $C_2$ as in Section 2.6 with parameter $\lambda_2$, there exists an $\varepsilon > 0$ such that

$$\mu_{1+\delta, \lambda_1(1+\delta)} \preceq \mu_{\lambda_2}^{(-, \varepsilon)}.$$

Hence, $(\mu_{\lambda_1}, \mu_{\lambda_2})$ is downward movable.    $\square$

For the rest of this section we will only consider the graph $\mathbb{T}^d$ for $d \geq 2$. The following is a 0–1 law for the upper invariant measure for the contact process.

PROPOSITION 4.2. *Let $\mathcal{A} \subseteq \{0,1\}^{\mathbb{T}^d}$, where $d \geq 2$, be a set which is invariant under all graph automorphisms on $\mathbb{T}^d$. Then, for $\lambda > 0$, we have that*

$$\mu_\lambda(\mathcal{A}) \in \{0,1\}.$$

PROOF. Let $\varepsilon > 0$. By elementary measure theory, there exists a cylinder event $\mathcal{B}$ depending on finitely many coordinates such that

$$(27) \qquad \mu_\lambda(\mathcal{A} \Delta \mathcal{B}) \leq \varepsilon.$$

Let $\operatorname{supp} \mathcal{B}$ denote the finite number of coordinates with respect to which $\mathcal{B}$ is measurable. Letting $\{T_\lambda(t)\}_{t \geq 0}$ denote the Markov semigroup for the contact process with parameter $\lambda$, we have that $\delta_1 T_\lambda(t) \to \mu_\lambda$ and also that $\mu_\lambda \preceq \delta_1 T_\lambda(t)$ for every $t \geq 0$. Choose $t$ so that, for all (equivalently, some) sites $s$,

$$\delta_1 T_\lambda(t)(\eta(s)=1) \leq \mu_\lambda(\eta(s)=1) + \frac{\varepsilon}{2|\operatorname{supp}\mathcal{B}|}.$$

It follows easily that if $m$ is any coupling of $\delta_1 T_\lambda(t)$ and $\mu_\lambda$ which is concentrated on $\{(\eta,\delta) : \eta \preceq \delta\}$, then, for any finite set $S$ of sites,

$$m((\eta,\delta) : \eta(s) \neq \delta(s) \text{ occurs for some } s \in S) \leq \frac{|S|\varepsilon}{2|\operatorname{supp}\mathcal{B}|}.$$

In particular, if $E$ is any event depending on at most $2|\operatorname{supp}\mathcal{B}|$ sites, then

$$(28) \qquad |\delta_1 T_\lambda(t)(E) - \mu_\lambda(E)| \leq \varepsilon.$$

For this fixed $t$, Theorem 4.6, page 35 of [19] shows that there exists an automorphism $\gamma \in AUT(\mathbb{T}^d)$ such that

$$(29) \qquad |\delta_1 T_\lambda(t)(\mathcal{B} \cap \gamma\mathcal{B}) - \delta_1 T_\lambda(t)(\mathcal{B})\delta_1 T_\lambda(t)(\gamma\mathcal{B})| \leq \varepsilon.$$

Furthermore, since $\mu_\lambda$ is invariant under automorphisms, (27) implies that

$$\mu_\lambda(\gamma\mathcal{A} \Delta \gamma\mathcal{B}) \leq \varepsilon,$$

and since $\mathcal{A} = \gamma\mathcal{A}$, we have

$$\mu_\lambda(\mathcal{A} \Delta \gamma\mathcal{B}) \leq \varepsilon.$$

It follows that

$$\mu_\lambda(\mathcal{B} \Delta \gamma\mathcal{B}) \leq \mu_\lambda(\mathcal{A} \Delta \gamma\mathcal{B}) + \mu_\lambda(\mathcal{A} \Delta \mathcal{B}) \leq 2\varepsilon.$$

Next, (28) implies that

$$|\delta_1 T_\lambda(t)(\mathcal{B}\Delta\gamma\mathcal{B}) - \mu_\lambda(\mathcal{B}\Delta\gamma\mathcal{B})| \le \varepsilon,$$

and so

(30) $$\delta_1 T_\lambda(t)(\mathcal{B}\Delta\gamma\mathcal{B}) \le 3\varepsilon.$$

We get that

$$\begin{aligned}
|\mu_\lambda(\mathcal{A}) - \mu_\lambda(\mathcal{A})^2| &= |\mu_\lambda(\mathcal{A}) - \mu_\lambda(\mathcal{A})\mu_\lambda(\gamma\mathcal{A})| \\
&\le |\mu_\lambda(\mathcal{B}) - \mu_\lambda(\mathcal{B})\mu_\lambda(\gamma\mathcal{B})| + 3\varepsilon \\
&\le |\delta_1 T_\lambda(t)(\mathcal{B}) - \delta_1 T_\lambda(t)(\mathcal{B})\delta_1 T_\lambda(t)(\gamma\mathcal{B})| + 6\varepsilon \\
&\le |\delta_1 T_\lambda(t)(\mathcal{B}) - \delta_1 T_\lambda(t)(\mathcal{B}\cap\gamma\mathcal{B})| + 7\varepsilon \\
&\le \delta_1 T_\lambda(t)(\mathcal{B}\Delta\gamma\mathcal{B}) + 7\varepsilon \le 10\varepsilon,
\end{aligned}$$

where we used (27), (28) and (29) for the three first inequalities and (30) in the last. Since $\varepsilon > 0$ was choosen arbitrarily, we get that

$$\mu_\lambda(\mathcal{A}) = \mu_\lambda(\mathcal{A})^2$$

and so $\mu_\lambda(\mathcal{A}) \in \{0, 1\}$.  □

REMARKS. The above proof works for any transitive and even quasi-transitive graph. For the case of $\mathbb{Z}^d$, this was proved in Proposition 2.16, page 143 of [19]. It is mentioned there that, while $\delta_1 T_\lambda(t)$ is ergodic for each $t$, one cannot conclude immediately the ergodicity of $\mu_\lambda$ because the class of ergodic processes is not weakly closed. We point out, however, that there is another important notion of convergence given by the $\bar{d}$-metric (see [24], page 89 for definition) on stationary processes. Convergence in this metric is stronger than weak convergence and weaker than convergence in the total variation norm. It is also known that the ergodic processes are $\bar{d}$-closed and that weak convergence together with stochastic ordering implies $\bar{d}$-convergence. In this way, one can conclude ergodicity of $\mu_\lambda$ using the $\bar{d}$-metric, giving an alternative proof of Proposition 2.16 of [19]. In fact, the proof of Proposition 4.2 is essentially based on this idea. However, because of the open question listed below, it is not so easy to formulate the $\bar{d}$-metric for tree indexed processes and so we choose a more hands on approach. Observe that the crucial property of $\bar{d}$-convergence which is essentially used in the above proof is that, for each fixed $k$, one has uniform convergence of the probability measures (in, say, the total variation norm) over all sets which depend on at most $k$ points. (The point is that the $k$ points can lie anywhere and, hence, this is much stronger than weak convergence.)

*Open question related to defining the $\bar{d}$-metric for tree indexed processes.* Assume that $\mu$ and $\nu$ are two automorphism invariant probability measures on $\{0,1\}^{\mathbb{T}^d}$ such that $\mu \preceq \nu$. Does there exist a $\mathbb{T}^d$-invariant coupling $(X, Y)$ with $X \sim \mu$, $Y \sim \nu$ and $X \preceq Y$?

PROPOSITION 4.3. *On $\mathbb{T}^d$, $d \geq 2$, there exists a $\lambda_p$ such that, for all $\lambda > \lambda_p$,*

$$\mu_\lambda(\mathcal{C}^+) = 1.$$

PROOF. By Theorem 1.33(c), page 275 in [19], for sufficiently large $\lambda$, $\mu_\lambda(\eta(s) = 1) \geq 2/3$. By [12], we have that if $\mu_\lambda(\eta(s) = 1) \geq 2/3$, then

$$\mu_\lambda(\mathcal{C}^+) > 0.$$

Finally, Proposition 4.2 then implies that

$$\mu_\lambda(\mathcal{C}^+) = 1. \qquad \square$$

**5. Relationship between $\varepsilon$-movability and dynamics.** In the general setup we have a family of stationary Markov processes parametrized by one or two parameters, for example, the contact processes $\Psi^\lambda$ ($\lambda$ is here the only parameter) or a stochastic Ising model $\Psi^{+,\beta,h}$ ($\beta$ and $h$ being the parameters). Many of the proofs in this paper will involve comparing the marginal distributions of these Markov processes for two different values of one of the involved parameters. Let $p$ be the parameter and let $p_1 < p_2$. Assume that the marginal distributions are $\mu_{p_1}$ and $\mu_{p_2}$, respectively, and that $\mu_{p_1} \preceq \mu_{p_2}$. Lemmas 5.1 and 5.2 show that there is a close connection between showing that $(\mu_{p_1}, \mu_{p_2})$ is downward $\varepsilon$-movable and that the infimum of the second process over a short time interval is stochastically larger than the first process.

Let $\Psi^\mu$ be a stationary Markov process on $\{0,1\}^S$ with marginal distribution $\mu$ and let $\{X_t\}_{t \geq 0} \sim \Psi^\mu$. For $\delta > 0$ and $s \in S$, define

$$X_{\inf,\delta}(s) := \inf_{t \in [0,\delta]} X_t(s),$$

and denote the distribution of $X_{\inf,\delta}$ by $\mu_{\inf,\delta}$. Similarly, define

$$X_{\sup,\delta}(s) := \sup_{t \in [0,\delta]} X_t(s),$$

and denote the distribution of $X_{\sup,\delta}$ by $\mu_{\sup,\delta}$.

LEMMA 5.1. *Take $S$ to be the sites of a bounded degree graph. Let $\{C(s,\sigma)\}_{s \in S, \sigma \in \{-1,1\}^S}$ be the flip rate intensities for a stationary Markov process $\Psi^\mu$ on $\{-1,1\}^S$ with marginal distribution $\mu$. Let*

$$\lambda := \sup_{(s,\sigma)} C(s,\sigma).$$

*For any $\tau > 0$, if we set $\varepsilon := 1 - e^{-\lambda\tau}$, we have that*

$$\mu^{(-,\varepsilon)} \preceq \mu_{\inf,\tau}.$$

*Similarly, we get that*

$$\mu_{\sup,\tau} \preceq \mu^{(+,\varepsilon)}.$$

PROOF. We will prove the first statement, the second statement follows by symmetry. Take $\tau > 0$. For every $s \in S$, associate an independent Poisson process with parameter $\lambda$. Define $\{(X_t^1, X_t^2)\}_{t \geq 0}$ in the following way. Let $X_0^1 \equiv X_0^2 \sim \mu$, and take $t'$ to be an arrival time for the Poisson process of a site $s$. For $i \in \{1, 2\}$, let $X_{t',-}^i$ and $X_{t',+}^i$ denote the configurations before and after the arrival. We let $X_{t',+}^1(s) \neq X_{t',-}^1(s)$ with probability $C(s, X_{t',-}^1)/\lambda$ and we let $X_{t',+}^2(s) = 0$ and finally, we let $X_{t',+}^1(S \setminus s) \equiv X_{t',-}^1(S \setminus s)$, $X_{t',+}^2(S \setminus s) \equiv X_{t',-}^2(S \setminus s)$. Do this independently for all arrival times for all Poisson processes of all sites. Observe that once $X_t^2(s)$ is 0, it remains so. Note also that $X_\tau^1 \sim \mu$, $X_\tau^2 \sim \mu^{(-,\varepsilon)}$. Furthermore, if $X_t^1(s) = 0$ for some $t \in [0, \tau]$, the construction guarantees that $X_\tau^2(s) = 0$ and, therefore, $X_\tau^2 \preceq X_{\inf,\tau}^1 \sim \mu_{\inf,\tau}$. □

LEMMA 5.2. *Take $S$ to be the sites of any bounded degree graph. Let $\{C(s, \sigma)\}_{s \in S, \sigma \in \{-1,1\}^S}$ be the flip rate intensities of a stationary Markov process $\Psi^\mu$ on $\{-1, 1\}^S$ with marginal distribution $\mu$. Define*

$$\lambda_1 := \inf_{s, \sigma : \sigma(s) = 1} C(s, \sigma).$$

*If $\lambda_1 > 0$, then for any $0 < \varepsilon < 1$, if we set $\tau := -\frac{\log(1-\varepsilon)}{\lambda_1}$, we have that*

$$\mu_{\inf,\tau} \preceq \mu^{(-,\varepsilon)}.$$

*Similarly, defining $\lambda_2 := \inf_{s, \sigma : \sigma(s) = 0} C(s, \sigma)$, if $\lambda_2 > 0$, then for any $0 < \varepsilon < 1$, if we set $\tau := -\frac{\log(1-\varepsilon)}{\lambda_2}$, we have that*

$$\mu^{(+,\varepsilon)} \preceq \mu_{\sup,\tau}.$$

PROOF. We will prove the first statement, the second statement follows by symmetry. For every $s \in S$, associate an independent Poisson process with parameter $\lambda := \sup_{(s,\sigma)} C(s, \sigma)$. Next, let $\{U_{s,k}\}_{s \in S, k \geq 1}$ be independent uniform $[0, 1]$ random variables also independent of the Poisson processes. If $t'$ is an arrival time for the Poisson process at site $s$, we write $U_{s,t'}$ for $U_{s,k}$, where $k$ is such that $t'$ is the time of the $k$th arrival of the Poisson process at site $s$. Define $\{(X_t^1, X_t^2)\}_{t \geq 0}$ in the following way. Let $X_0^1 \equiv X_0^2 \sim \mu$, and take $t'$ to be an arrival time for the Poisson process of a site $s$. We let $X_{t',+}^1(s) \neq X_{t',-}^1(s)$ if $U_{s,t'} \leq C(s, X_{t',-}^1)/\lambda$. Furthermore, we let $X_{t',+}^2(s) = 0$ if $U_{s,t'} \leq \lambda_1/\lambda$ or $X^2_{t',-}(s) = 0$, and finally, we let $X^1_{t',+}(S \setminus s) \equiv X^1_{t',-}(S \setminus s)$, $X^2_{t',+}(S \setminus s) \equiv X^2_{t',-}(S \setminus s)$. Do this independently for all arrival times for all Poisson processes of all sites. Clearly, $X^1_\tau \sim \mu$ and $X^2_\tau \sim \mu^{(-,\varepsilon)}$. Furthermore, if $X^2_\tau(s) = 0$, then either $X^1_0(s) = X^2_0(s) = 0$ or there exists a $t \in [0,\tau]$ such that $t$ is an arrival time for the Poisson process associated to $s$ and $U_{s,t} \leq \lambda_1/\lambda$. Since $\lambda_1 \leq C(s, X^1_{t^-})$ if $X^1_{t^-}(s) = 1$, we get that either $X^1_{t^+}(s)$ or $X^1_{t^-}(s)$ is 0 and, therefore, $X^1_{\inf,\tau} \preceq X^2_\tau$. $\square$

To illustrate why the condition $\lambda_1 > 0$ of Lemma 5.2 is needed, consider the case $\mu = \pi_p$ for some $p > 0$. With $\varepsilon > 0$, if we assume the trivial dynamics $C(s, \sigma) = 0$ for all $s, \sigma$, we will of course not have that $\mu_{\inf,\tau} \preceq \mu^{(-,\varepsilon)}$ for any $\tau > 0$.

**6. Proof of Theorem 1.9.** Take $\lambda > \lambda_p$ and let $\lambda' = (\lambda + \lambda_p)/2$. By Theorem 1.11, there exists an $\varepsilon > 0$ such that $(\mu_{\lambda'}, \mu_\lambda)$ is downward $\varepsilon$-movable. Lemma 5.1 gives us that there exists a $\tau > 0$ such that $\mu_{\lambda'}^{(-,\varepsilon)} \preceq \mu_{\lambda,\inf,\tau}$ and, hence, that $\mu_{\lambda'} \preceq \mu_{\lambda,\inf,\tau}$. Therefore, since $\mathcal{C}^+$ is an increasing event and $\lambda' > \lambda_p$, we have that

$$1 = \mu_{\lambda'}(\mathcal{C}^+) \leq \mu_{\lambda,\inf,\tau}(\mathcal{C}^+)$$

and so

$$\Psi^\lambda(\mathcal{C}^+_t \ \forall t \in [0,\tau]) = 1.$$

The theorem now follows from countable additivity. $\square$

**7. Proof of Theorem 1.1.** In this section we will deal with stationary distributions for interacting particle systems which are monotone in the sense of Definition 2.2.

Let $G = (S, E)$ be a countable connected locally finite graph and let $\Lambda \subseteq S$ be connected and $|\Lambda| < \infty$. Let $\{\mu^p_\Lambda\}_{p \in I}$, where $I \subseteq \mathbb{R}$ be a family of probability measures on $\{-1, 1\}^\Lambda$ such that

$$\mu^{p_1}_\Lambda \preceq \mu^{p_2}_\Lambda \qquad \forall p_1 \leq p_2.$$

Assume that there exist stationary Markov processes $\Psi^p_\Lambda$ governed by flip rate intensities $\{C_{p,\Lambda}(s, \sigma)\}_{s \in \Lambda, \sigma \in \{-1,1\}^\Lambda}$ and with marginal distributions $\mu^p_\Lambda$. Furthermore, assume that there exists limiting distributions $\Psi^p$ of $\Psi^p_\Lambda$ and $\mu^p$ of $\mu^p_\Lambda$ as $\Lambda \uparrow S$. Assume that $\mu^p_\Lambda$ are monotone for every $p$ and $\Lambda$. For $p_1 < p_2$, let

$$A_{\Lambda,p_1,p_2} := \inf_{\substack{s \in \Lambda \\ \xi \in \{-1,1\}^{\Lambda \setminus s}}} [\mu^{p_2}_\Lambda(\sigma(s) = 1 | \sigma(\Lambda \setminus s) \equiv \xi) - \mu^{p_1}_\Lambda(\sigma(s) = 1 | \sigma(\Lambda \setminus s) \equiv \xi)]$$

and assume that, for all $p_1 < p_2$,

$$\inf_{\Lambda \subseteq S} A_{\Lambda, p_1, p_2} > 0.$$

For fixed $p_1 < p_2$, there exists by Proposition 3.4 an $\varepsilon > 0$ such that $(\mu^{p_1}, \mu^{p_2})$ is both upward and downward $\varepsilon$-movable. Next, by Lemma 5.1, there exists a $\tau > 0$ such that

$$\mu^{p_2, (-, \varepsilon)} \preceq \mu^{p_2}_{\inf, \tau},$$

and therefore,

$$\mu^{p_1} \preceq \mu^{p_2}_{\inf, \tau}. \tag{31}$$

THEOREM 7.1. *Consider the setup just described. Let $\mathcal{A}$ be an increasing event on $\{-1, 1\}^S$ and let $\mathcal{A}_t$ be the event that $\mathcal{A}$ occurs at time $t$.*

(1) *Let $a \in \mathbb{R}$. If*

$$\mu^p(\mathcal{A}) = 1$$

*for all $p \in I$ with $p > a$, then*

$$\Psi^p(\mathcal{A}_t \text{ occurs for every } t) = 1$$

*for all $p \in I$ with $p > a$.*

(2) *Let $a \in \mathbb{R}$. If*

$$\mu^p(\mathcal{A}) = 0$$

*for all $p \in I$ with $p < a$, then*

$$\Psi^p(\mathcal{A}_t \text{ occurs for some } t) = 0$$

*for all $p \in I$ with $p < a$.*

PROOF. We prove only (1), as (2) is proved in an identical way. Take $p > a$ and let $p_2 = (p + a)/2$. By the argument leading toward (31), there exists $\tau > 0$ such that

$$\mu^{p_2}(\mathcal{A}) \le \mu^p_{\inf, \tau}(\mathcal{A}).$$

By using $\mu^{p_2}(\mathcal{A}) = 1$ and

$$\mu^p_{\inf, \tau}(\mathcal{A}) \le \Psi^p(\mathcal{A}_t \text{ occurs for every } t \in [0, \tau]),$$

we get by countable additivity that

$$\Psi^p(\mathcal{A}_t \text{ occurs for every } t) = 1. \qquad \square$$

We will now be able to prove Theorem 1.1 easily.

PROOF OF THEOREM 1.1. We prove only the very first statement; all the other statements are proved in a similar manner. We fix $\beta \geq 0$ and then $h$ will correspond to our parameter $p$ in the above set up. For any $\Lambda \subseteq S$, any $s \in \Lambda$ and any $\xi \in \{-1, 1\}^{\Lambda \setminus s}$, we have that

$$(32) \qquad \mu_{\Lambda}^{+, \beta, h}(\sigma(s) = 1 | \sigma(\Lambda \setminus s) = \xi) = \frac{1}{1 + e^{-2\beta(\sum_{t \,:\, t \sim s} \xi(t)) - 2h}},$$

where we let $\xi(t) = 1$ if $t \in \Lambda^c$ in order to take the boundary condition into account. It is obvious from (32) and the definition of monotonicity that $\mu_{\Lambda}^{+, \beta, h}$ is monotone for any $h$ and $\Lambda$. Letting $h_1 < h_2$, it is immediate that

$$A_{\Lambda, h_1, h_2} = \inf_{\substack{s \in \Lambda \\ \xi \in \{-1, 1\}^{\Lambda \setminus s}}} \left[ \frac{1}{1 + e^{-2\beta(\sum_{t \,:\, t \sim s} \xi(t)) - 2h_2}} - \frac{1}{1 + e^{-2\beta(\sum_{t \,:\, t \sim s} \xi(t)) - 2h_1}} \right] > 0,$$

where again $\xi(t) = 1$ for all $t \in \Lambda^c$. It is not hard to see that this strict inequality must hold uniformly in $\Lambda$, that is,

$$\inf_{\Lambda \subseteq S} A_{\Lambda, h_1, h_2} > 0.$$

It follows that all of the assumptions of Theorem 7.1 hold and part (1) of that result gives us what we want. $\square$

PROOF OF LEMMA 1.2. Fix $\beta \geq 0$. Given any $p \in (0, 1)$, it is easy to see that there exists a real number $h_2$ such that, for all $h \geq h_2$, for $s \in S$ and for all $\xi \in \{-1, 1\}^{S \setminus s}$,

$$\mu^{+, \beta, h}(\sigma(s) = 1 | \sigma(S \setminus s) = \xi) \geq p$$

and, hence, $\pi_p \preceq \mu^{+, \beta, h}$. It is also easy to see that there exists a real number $h_1$ such that, for all $h < h_1$, for $s \in S$ and for all $\xi \in \{-1, 1\}^{S \setminus s}$,

$$\mu^{+, \beta, h}(\sigma(s) = 1 | \sigma(S \setminus s) = \xi) \leq p$$

and, hence, $\mu^{+, \beta, h} \preceq \pi_p$. The statements of the lemma easily follow from these facts. $\square$

## 8. Proof of Theorem 1.3.
In this section we will use a variant of the so-called Peierls argument to prove Theorem 1.3. We prove this only for $\mathbb{Z}^2$; the proof (with more complicated topological details) can be carried out for $\mathbb{Z}^d$ with $d \geq 3$.

We will write $0 \overset{-, t}{\longleftrightarrow} \partial \Lambda_L$ for the event that there exists a path of sites in state $-1$ connecting the origin to $\partial \Lambda_L := \Lambda_{L+1} \setminus \Lambda_L$ at time $t$ and we will write $0 \overset{-, t}{\longleftrightarrow} \infty$ for the event that there exists an infinite path of sites in state $-1$ containing the origin at time $t$. We will also write $0 \overset{+,t}{\longleftrightarrow} \partial \Lambda_L$ and $0 \overset{+,t}{\longleftrightarrow} \infty$ for the obvious analogous events. We will first need Lemma 8.1 and the concept of a dual graph. The dual graph $G_n^{\mathrm{dual}} = (S_n^{\mathrm{dual}}, E_n^{\mathrm{dual}})$ of $G_n = (S_n, E_n)$ consists of the set of sites $S_n^{\mathrm{dual}} := \{-n - \frac{1}{2}, \ldots, n + \frac{1}{2}\}^2$ and $E_n^{\mathrm{dual}}$, which is the set of nearest neighbor pairs of $S_n^{\mathrm{dual}}$. In this paper we will only work with the edges of the dual graph. An edge $e \in E_n^{\mathrm{dual}}$ crosses one (and only one) edge $f \in E_n$ and the end sites of this edge $f$ will be called the sites (of $G_n$) associated to $e$. For a random spin configuration $X$ on $\{-1, 1\}^{S_n}$, define a random edge configuration $Y$ on $\{0, 1\}^{E_n^{\mathrm{dual}}}$ in the following way:

$$(33) \qquad Y(e) = \begin{cases} 0, & \text{if } X(t) = X(s), \\ 1, & \text{if } X(t) \neq X(s), \end{cases}$$

where $s, t$ are the sites associated to edge $e \in E_n^{\mathrm{dual}}$. In Figure 1 we have drawn a configuration $\sigma \in \{-1, 1\}^{S_1}$ and the induced edge configuration on $\{0, 1\}^{E_1^{\mathrm{dual}}}$.

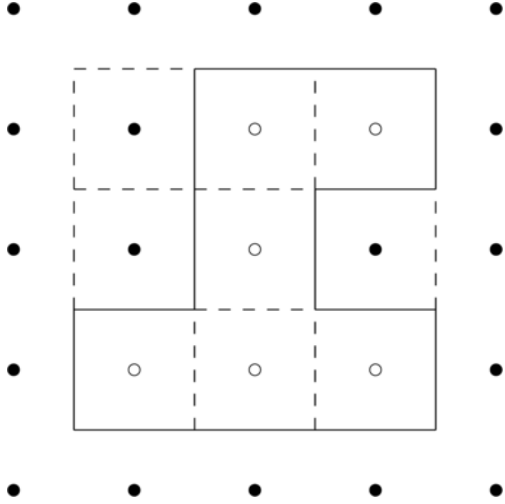

Fig. 1. $S_1$ and the edges of its dual graph. A solid circle marks a site with spin 1, while an empty circle has spin $-1$. A solid line is a present edge of the dual graph, and a dashed line is an absent edge of the dual graph.

Assume that the sites evolve according to the flip rate intensities $\{C_n(s, \sigma)\}_{s \in S_n, \sigma \in \{-1, 1\}^{S_n}}$. Consider $\gamma$, a (finite) path of edges in the dual graph. Take $\gamma'$ to be a subset of $\gamma$. Assume that all edges of $\gamma'$ are absent and all edges of $\gamma \setminus \gamma'$ are present at $t = 0$. We want to estimate the probability of the event that all edges of $\gamma'$ are present at some point (not necessarily all at the same time) during some time interval $[0, \tau]$. In other words, we want to estimate $\mathbb{P}(Y_{\sup, \tau}(\gamma') \equiv 1 | Y_0(\gamma') \equiv 0, Y_0(\gamma \setminus \gamma') \equiv 1)$.

LEMMA 8.1. *Let* $\{C_n(s,\sigma)\}_{s\in S_n, \sigma\in\{-1,1\}^{S_n}}$ *be the flip rate intensities for a stationary Markov process on* $\{-1,1\}^{S_n}$ *and let* $Y_t$ *be defined as above. Let*

$$\lambda := \sup_{(s,\sigma)} C_n(s,\sigma) \ (<\infty).$$

*For any* $\tau > 0$ *and any* $\gamma' \subseteq E_n^{\text{dual}}$,

$$\mathbf{P}(Y_{\sup,\tau}(\gamma') \equiv 1 | Y_0(\gamma') \equiv 0, Y_0(E_n^{\text{dual}} \setminus \gamma') \equiv 1) \le (4(1 - e^{-\lambda\tau})^{1/4})^{|\gamma'|}.$$

PROOF. Take $\tau > 0$. For every $s \in S_n$, associate an independent Poisson process with parameter $\lambda$. Define $\{X_t\}_{t\geq 0}$ in the following way. Let $X_0 \sim \mu$ and take $t'$ to be an arrival time for the Poisson process of a site $s$. We let $X_{t',+}(s) \neq X_{t',-}(s)$ with probability $C(s, X_{t',-})/\lambda$. Do this independently for all arrival times for all Poisson processes associated to the different sites. It is immediate that $X_\tau \sim \mu$. Let $s_i, i \in \{1,\dots,l\}$, be distinct sites of $S_n$. The event $\{X_{\inf,\tau}(s_i) \neq X_{\sup,\tau}(s_i) \ \forall i \in \{1,\dots,l\}\}$ is contained in the event that every Poisson process associated to the sites $s_i, i \in \{1,\dots,l\}$, have had at least one arrival by time $\tau$. The probability that a particular site has had an arrival by time $\tau$ is $1 - e^{-\lambda\tau}$. Furthermore, this event is independent of the Poisson processes for all other sites. Therefore,

$$(34) \qquad \mathbf{P}(X_{\inf,\tau}(s_i) \neq X_{\sup,\tau}(s_i) \ \forall i \in \{1,\dots,l\}) \le (1 - e^{-\lambda\tau})^l.$$

Given $\gamma'$, consider the set of all sites associated to some edge of $\gamma'$ and let $n_{\gamma'}$ be the cardinality of that set. Observe that $n_{\gamma'} \le 2|\gamma'|$ and that in order for the event $(Y_{\sup,\tau}(\gamma') \equiv 1 | Y_0(\gamma') \equiv 0, Y_0(E_n^{\text{dual}} \setminus \gamma') \equiv 1)$ to occur, at least $|\gamma'|/4$ of the sites associated to $\gamma'$ must flip during $[0,\tau]$. This is because one site is associated to at most 4 edges. Denote the event that at least $|\gamma'|/4$ of the sites associated to $\gamma'$ flip during $[0,\tau]$ by $\mathcal{A}_{\tau,\gamma'}$. Take $\tilde{S}$ to be a subset of the sites associated to $\gamma'$ such that $|\tilde{S}| \ge |\gamma'|/4$. By (34), the probability that all of these sites flip during $[0,\tau]$ is less than $(1 - e^{-\lambda\tau})^{|\tilde{S}|} \le (1 - e^{-\lambda\tau})^{|\gamma'|/4}$. To conclude, observe that the number of subsets of the sites associated to $\gamma'$ is bounded by $2^{2|\gamma'|}$. Hence, the probability of the event $\mathcal{A}_{\tau,\gamma'}$ must be less than $(1 - e^{-\lambda\tau})^{|\gamma'|/4} 2^{2|\gamma'|}$, and so

$$\mathbf{P}(Y_{\sup,\tau}(\gamma') \equiv 1 | Y_0(\gamma') \equiv 0, Y_0(E_n^{\text{dual}} \setminus \gamma') \equiv 1)$$
$$\le \mathbf{P}(\mathcal{A}_{\tau,\gamma'}) \le ((1 - e^{-\lambda\tau})^{1/4}4)^{|\gamma'|}. \qquad \square$$

PROOF OF THEOREM 1.3. We will prove the theorem for $d = 2$. For $\beta > \beta_p$, choose $\delta_1 > 0$ so that $\beta' := \beta\frac{2-\delta_1}{2} > \beta_p$ and, hence,

$$\sum_{l=1}^{\infty} l 3^{l-1} e^{-2\beta' l} < \infty.$$

Next, choose $N$ and $\varepsilon < 1/2$ such that $\frac{4}{N} \le \delta_1$, and $\varepsilon^{1/N} \le e^{-\beta(2-\delta_1)}$ and let $\tau$ be such that $\varepsilon = 4(1 - e^{-\lambda\tau})^{1/4}$. Let $\delta > 0$ be arbitrary and choose $L$ so that

$$3\sum_{l=L}^{\infty} l 3^{l-1} e^{-2\beta' l} < \delta.$$

Let $\mathcal{E}_{L,\tau}$ be the event that $0 \xleftrightarrow{-,t} \partial\Lambda_L$, for some $t \in [0,\tau]$. Let $\Psi_n^{+,\beta}$ be defined as in Section 2.3. We will show that

$$\Psi_n^{+,\beta}(\mathcal{E}_{L,\tau}) < \delta \qquad \forall n > L.$$

Since $\Psi_n^{+,\beta}(\mathcal{E}_{L,\tau}) \to \Psi^{+,\beta}(\mathcal{E}_{L,\tau})$ (see Section 2.3) we get that $\Psi^{+,\beta}(\mathcal{E}_{L,\tau}) \le \delta$. Letting $L \to \infty$ and $\delta \to 0$, we get that

$$\Psi^{+,\beta}(\exists t \in [0,\tau]: 0 \xleftrightarrow{-,t} \infty) = 0,$$

and then by countable additivity,

$$\Psi^{+,\beta}(\exists t \ge 0: 0 \xleftrightarrow{-,t} \infty) = 0.$$

It is well known (see [8]) that if all sites in $\Lambda_{n+1} \setminus \Lambda_n$ take the value $+1$,

$$(35) \qquad \begin{aligned} \mathcal{E}_{L,\tau} &\subseteq \{\exists \gamma \subseteq E_n^{\text{dual}}, t \in [0,\tau] : |\gamma| \ge L, \gamma \text{ surrounds the origin}, Y_t(\gamma) \equiv 1\} \\ &\subseteq \{\exists \gamma \subseteq E_n^{\text{dual}} : |\gamma| \ge L, \gamma \text{ surrounds the origin}, Y_{\sup,\tau}(\gamma) \equiv 1\}. \end{aligned}$$

To prove $\Psi_n^{+,\beta}(\mathcal{E}_{L,\tau}) < \delta$, consider $\gamma$ with $|\gamma| = l$ a contour in $E_n^{\text{dual}}$ surrounding the origin. By Lemma 8.1, $\mathbf{P}(Y_{\sup,\tau}(\gamma') \equiv 1 | Y_0(\gamma') \equiv 0, Y_0(\gamma \setminus \gamma') \equiv 1) \le \varepsilon^{|\gamma'|}$ whenever $\gamma' \subseteq \gamma$. We get

$$\begin{aligned} &\mathbf{P}(Y_{\sup,\tau}(\gamma) \equiv 1) \\ &= \sum_{k=0}^{l} \sum_{\substack{\gamma' \subseteq \gamma \\ |\gamma'|=k}} \mathbf{P}(Y_0(\gamma') \equiv 0, Y_0(\gamma \setminus \gamma') \equiv 1) \\ &\qquad\qquad \times \mathbf{P}(Y_{\sup,\tau}(\gamma') \equiv 1 | Y_0(\gamma') \equiv 0, Y_0(\gamma \setminus \gamma') \equiv 1) \\ (36) \qquad &\le \sum_{k=0}^{l} \sum_{\substack{\gamma' \subseteq \gamma \\ |\gamma'|=k}} \mathbf{P}(Y_0(\gamma') \equiv 0, Y_0(\gamma \setminus \gamma') \equiv 1)\varepsilon^k \\ &= \sum_{k=0}^{l} \mathbf{P}(\{\text{all edges except } k \text{ of } \gamma \text{ are present at } t=0\})\varepsilon^k \\ &= \sum_{k=0}^{l/N} \mathbf{P}(\{\text{all edges except } k \text{ of } \gamma \text{ are present at } t=0\})\varepsilon^k \end{aligned}$$

$$+ \sum_{k=l/N+1}^{l} \mathbf{P}(\{\text{all edges except } k \text{ of } \gamma \text{ are present at } t=0\})\varepsilon^k.$$

Obviously, $l/N$ need not be an integer, but correcting for this is trivial and is left for the reader.

We need to estimate $\mathbf{P}(\{\text{all edges except } k \text{ of } \gamma \text{ are present at } t=0\})$. For this purpose, define $T: \{-1,1\}^{S_n} \to \{-1,1\}^{S_n}$, by

$$(T\sigma)(s) = \begin{cases} \sigma(s), & \text{if } s \text{ is not in the domain bounded by } \gamma, \\ -\sigma(s), & \text{if } s \text{ is in the domain bounded by } \gamma, \end{cases}$$

for all $\sigma \in \{-1,1\}^{S_n}$. Let $E_k = \{\sigma: \text{all edges except } k \text{ of } \gamma \text{ are present}\}$. Since $H_n^{+,\beta}$ of (6) gives a contribution of $-\beta$ for adjacent pairs of equal spin and $+\beta$ for adjacent pairs of unequal spin, we have that, for $\sigma \in E_k$, $H_n^{+,\beta}(T\sigma) = H_n^{+,\beta}(\sigma) - 2\beta(|\gamma| - k) + 2\beta k = H_n^{+,\beta}(\sigma) - 2\beta|\gamma| + 4\beta k$.

Hence, for $\sigma \in E_k$,

$$\mu_n^{+,\beta}(\sigma) = \frac{e^{-H_n^{+,\beta}(\sigma)}}{Z} = \frac{e^{-H_n^{+,\beta}(T\sigma) - 2\beta|\gamma| + 4\beta k}}{Z},$$

and so

$$\mu_n^{+,\beta}(E_k) = \sum_{\sigma \in E_k} \mu_n^{+,\beta}(\sigma) = e^{-2\beta l + 4\beta k} \sum_{\sigma \in E_k} \frac{e^{-H_n^{+,\beta}(T\sigma)}}{Z}$$

$$\leq e^{-2\beta l + 4\beta k} \sum_{\sigma \in \{-1,1\}^{S_n}} \frac{e^{-H_n^{+,\beta}(T\sigma)}}{Z} = e^{-2\beta l + 4\beta k},$$

where the last equality follows from $T$ being bijective. We then get that

$$\sum_{k=0}^{l/N} \mathbf{P}(\{\text{all edges except } k \text{ of } \gamma \text{ are present at } t=0\})\varepsilon^k$$

(37)
$$\leq \sum_{k=0}^{l/N} e^{-2\beta l + 4\beta k} \varepsilon^k \leq e^{-2\beta l + 4\beta l/N} \sum_{k=0}^{l/N} \varepsilon^k \leq 2 e^{-2\beta l + 4\beta l/N}$$

$$\leq 2 e^{-\beta(2-\delta_1)l} = 2 e^{-2\beta' l}.$$

Furthermore,

$$\sum_{k=l/N+1}^{l} \mathbf{P}(\{\text{all edges except } k \text{ of } \gamma \text{ are present at } t=0\})\varepsilon^k$$

(38)
$$\leq \varepsilon^{l/N} \sum_{k=l/N+1}^{l} \mathbf{P}(\{\text{all edges except } k \text{ of } \gamma \text{ are present at } t=0\})$$

$$\leq \varepsilon^{l/N} \leq e^{-\beta(2-\delta_1)l} = e^{-2\beta' l},$$

where we use that {all edges except $k$ of $\gamma$ are present at $t = 0$} are disjoint events for different $k$. Hence, (36), (37) and (38) combined give us

$$\mathbf{P}(Y_{\sup,\tau}(\gamma) \equiv 1) \leq 3e^{-2\beta' l}$$

and so by (35), for all $n > L$,

$$\Psi_n^{+,\beta}(\mathcal{E}_{L,\tau}) \leq \Psi_n^{+,\beta}(\exists\, \gamma \subseteq E_n^{\mathrm{dual}} : |\gamma| \geq L, \gamma \text{ surrounds the origin}, Y_{\sup,\tau}(\gamma) \equiv 1)$$

$$\leq \sum_{l=L}^{\infty} l 3^{l-1} 3e^{-2\beta' l} < \delta,$$

where the second to last inequality follows from the fact that the number of contours around the origin of length $l$ is at most $l3^{l-1}$ (see [8]).  □

REMARK. For $\mathbb{Z}^d$, the proof is generalized by noting that the number of connected surfaces of size $l$ surrounding the origin is at most $C(d)^l$, for some constant $C(d)$. The arguments are the same but the "topological details" are messier.

**9. Proof of Theorem 1.5.** We will start this subsection by presenting a theorem by Liggett, Schonmann and Stacey [21].

THEOREM 9.1. *Let $G = (S, E)$ be a graph with a countable set of sites in which every site has degree at most $\Delta \geq 1$, and in which every finite connected component of $G$ contains a site of degree strictly less than $\Delta$. Let $p, \alpha, r \in [0, 1]$, $q = 1 - p$, and suppose that*

$$(1 - \alpha)(1 - r)^{\Delta-1} \geq q,$$

$$(1 - \alpha)\alpha^{\Delta-1} \geq q.$$

*If $\mu \in G(p)$, then $\pi_{\alpha r} \preceq \mu$. In particular, if $q \leq (\Delta - 1)^{\Delta-1}/\Delta^{\Delta}$, then $\pi_\rho \preceq \mu$, where*

$$\rho = \left(1 - \frac{q^{1/\Delta}}{(\Delta - 1)^{(\Delta-1)/\Delta}}\right)(1 - (q(\Delta - 1))^{1/\Delta}).$$

Here $G(p)$ denotes the set of probability measures on $\{-1, 1\}^S$ such that if $\mu \in G(p)$, $X \sim \mu$, then for any site $s \in S$,

$$\mathbb{P}[X(s) = 1 | \sigma(\{X(t) : \{s, t\} \notin E\})] \geq p \qquad \text{a.s.}$$

Observe that when $p \to 1 \Rightarrow q \to 0$ and so $\rho \to 1$. The above theorem is stated as the original in [21]. However, by considering the line-graph of $G = (S, E)$, it can be restated in the following way.

COROLLARY 9.2. *Let $\tilde{G} = (\tilde{S}, \tilde{E})$ be any countable graph of degree at most $\Delta$. For each $0 < \rho < 1$, there exists a $0 < p < 1$, where $p = p(\Delta, \rho)$ such that if $Y \sim \nu$, where $\nu$ is a probability measure on the edges of $\tilde{G}$ such that for every edge $e \in \tilde{E}$,*

$$\mathbb{P}[Y(e) = 1 | \sigma(\{Y(f) : e \not\sim f\})] \geq p \qquad a.s.,$$

*we have that $\pi_\rho^{\tilde{E}} \preceq \nu$.*

By $e \not\sim f$, we of course mean that the edges $e$ and $f$ do not have any endpoints in common. Here, $\pi_\rho^{\tilde{E}}$ is the product measure with density $\rho$ on the edges of $\tilde{G}$.

Consider a graph $G = (S, E)$ and a subgraph $G' = (S', E')$, where $S' = S$ and $E' \subset E$. Let $X \sim \pi_p$ on $S$. We declare an edge $e \in E'$ to be closed if any of the endpoints takes the value 0 under $X$. Corollary 9.2 gives us that, for any $\rho < 1$, there is a $p < 1$ such that this method of closing edges dominates independent bond percolation with density $\rho$ on $E'$. Observe that we can choose $p$ independent of $E'$ since the maximal degree of $E'$ is bounded above by the maximal degree of $E$.

Let $(X, Y) \sim \mathbf{P}_n^p$, defined in Section 2.5. Close every $e \in E_n$ such that $Y(e) = 1$ independently with probability $\varepsilon$, thus creating $(X, Y^{(-,\varepsilon)})$. Compare this to closing every site in $S_n$ independently with parameter $\varepsilon'$ [creating $X^{(-,\varepsilon')}$] and defining

$$Y^{\varepsilon'}(e) = \begin{cases} 1, & \text{if } Y(e) = 1 \text{ and neither one of the endpoints of } e \text{ flips}, \\ 0, & \text{otherwise}. \end{cases}$$

By the arguments of the last paragraph, we see that, for a fixed $\varepsilon$, there exists an $\varepsilon'$ [that we can choose independent of $(X, Y)$ and $n$] such that the first way (i.e., independent bond percolation) of removing edges is stochastically dominated by the latter. Hence,

$$\mathbf{P}_n^p((X, Y^{(-,\varepsilon)}) \in (\{-1, 1\}^{S_n}, \cdot) | (X, Y))$$
$$\preceq \mathbf{P}_n^p((X^{(-,\varepsilon')}, Y^{\varepsilon'}) \in (\{-1, 1\}^{S_n}, \cdot) | (X, Y)).$$

By averaging over all possible $(X, Y)$, the next lemma follows.

LEMMA 9.3. *With notation as above, for any $\varepsilon > 0$, there exists $\varepsilon' > 0$ independent of $n$ such that*

$$\mathbf{P}_n^p((X, Y^{(-,\varepsilon)}) \in (\{-1, 1\}^{S_n}, \cdot)) \preceq \mathbf{P}_n^p((X^{(-,\varepsilon')}, Y^{\varepsilon'}) \in (\{-1, 1\}^{S_n}, \cdot)).$$

Observe that

$$(39) \qquad \mathbf{P}_n^p((X, Y^{(-,\varepsilon)}) \in (\{-1, 1\}^{S_n}, \cdot)) \overset{\mathcal{D}}{=} \tilde{\nu}_n^{p,(-,\varepsilon)}(\cdot)$$

and that

$$(40) \qquad \mathbf{P}_n^p((X^{(-,\varepsilon')}, Y^{\varepsilon'}) \in (\cdot, \{-1,1\}^{E_n})) \stackrel{\mathcal{D}}{=} \mu_n^{+,\beta,(-,\varepsilon')}(\cdot).$$

We are now ready to prove Theorem 1.5.

PROOF OF THEOREM 1.5. For any choice of $\beta > \beta_c$, take $p = 1 - e^{-2\beta}$ and let $\delta \in (0, p - p_c)$. Now, (14) and Holley's inequality imply that

$$\tilde{\nu}_n^{p-\delta} \preceq \tilde{\nu}_n^p \qquad \forall\, n \in \mathbb{N}^+.$$

Since, by (14), both $\tilde{\nu}_n^{p-\delta}$ and $\tilde{\nu}_n^p$ are monotone, there exists by Lemma 3.3 (it is easy to check that all other conditions of that lemma are satisfied) an $\varepsilon > 0$ such that

$$(41) \qquad \tilde{\nu}_n^{p-\delta} \preceq \tilde{\nu}_n^{p,(-,\varepsilon)} \qquad \forall\, n \in \mathbb{N}^+.$$

In [13] they show that the limit $\lim_n \tilde{\nu}_n^{p-\delta}(0 \longleftrightarrow \partial \Lambda_n)$ exists and that

$$(42) \qquad \lim_n \tilde{\nu}_n^{p-\delta}(0 \longleftrightarrow \partial \Lambda_n) > 0.$$

Here $\{0 \longleftrightarrow \partial \Lambda_n\}$ denotes the event that there exists a path of present edges connecting the origin to $\partial \Lambda_n := \Lambda_{n+1} \setminus \Lambda_n$. Since $\{0 \longleftrightarrow \partial \Lambda_n\}$ is an increasing event on the edges, Lemma 9.3 guarantees the existence of an $\varepsilon' > 0$ such that

$$\tilde{\nu}_n^{p,(-,\varepsilon)}(0 \longleftrightarrow \partial \Lambda_n)$$
$$= \mathbf{P}_n^p((X, Y^{(-,\varepsilon)}) \in (\{-1,1\}^{S_n}, 0 \longleftrightarrow \partial \Lambda_n))$$
$$\leq \mathbf{P}_n^p((X^{(-,\varepsilon')}, Y^{\varepsilon'}) \in (\{-1,1\}^{S_n}, 0 \longleftrightarrow \partial \Lambda_n)) \qquad \forall\, n \in \mathbb{N}^+.$$

If there exists a path of present edges connecting the origin to the boundary $\partial \Lambda_n$ under $Y$, all the sites of this path must have the value 1 under $X$. Similarly for $(X^{(-,\varepsilon')}, Y^{\varepsilon'})$, if there exists a path of present edges connecting the origin to the boundary $\partial \Lambda_n$ under $Y^{\varepsilon'}$, all the sites of this path must have the value 1 under $X^{(-,\varepsilon')}$. Hence,

$$\mathbf{P}_n^p((X^{(-,\varepsilon')}, Y^{\varepsilon'}) \in (\{-1,1\}^{S_n}, 0 \longleftrightarrow \partial \Lambda_n))$$
$$= \mathbf{P}_n^p((X^{(-,\varepsilon')}, Y^{\varepsilon'}) \in (0 \stackrel{+}{\longleftrightarrow} \partial \Lambda_n, 0 \longleftrightarrow \partial \Lambda_n))$$
$$\leq \mathbf{P}_n^p((X^{(-,\varepsilon')}, Y^{\varepsilon'}) \in (0 \stackrel{+}{\longleftrightarrow} \partial \Lambda_n, \{0,1\}^{E_n}))$$
$$= \mu_n^{+,\beta,(-,\varepsilon')}(0 \stackrel{+}{\longleftrightarrow} \partial \Lambda_n).$$

Of course,

$$\mu_n^{+,\beta,(-,\varepsilon')}(0 \stackrel{+}{\longleftrightarrow} \partial \Lambda_n) \leq \mu_n^{+,\beta,(-,\varepsilon')}(0 \stackrel{+}{\longleftrightarrow} \partial \Lambda_L) \qquad \forall\, L < n.$$

Therefore, for any $L$, we have that

$$0 < \lim_n \tilde{\nu}_n^{p-\delta}(0 \longleftrightarrow \partial\Lambda_n)$$

$$\leq \lim_n \mu_n^{+,\beta,(-,\varepsilon')}(0 \xleftrightarrow{+} \partial\Lambda_L) = \mu^{+,\beta,(-,\varepsilon')}(0 \xleftrightarrow{+} \partial\Lambda_L),$$

and so

$$0 < \lim_L \mu^{+,\beta,(-,\varepsilon')}(0 \xleftrightarrow{+} \partial\Lambda_L) = \mu^{+,\beta,(-,\varepsilon')}(0 \xleftrightarrow{+} \infty).$$

The limit in $L$ exists since $\{0 \xleftrightarrow{+} \partial\Lambda_{L_2}\} \subseteq \{0 \xleftrightarrow{+} \partial\Lambda_{L_1}\}$ for $L_1 \leq L_2$. Since $\mu^{+,\beta}$ is ergodic (see [19], pages 143 and 195), it follows that $\mu^{+,\beta,(-,\varepsilon')}$ must also be ergodic. This is because $\mu^{+,\beta,(-,\varepsilon')}$ can be expressed as a function of two independent processes, one being $\mu^{+,\beta}$ and the other a product measure. We conclude that

$$(43) \qquad \mu^{+,\beta,(-,\varepsilon')}(\mathcal{C}^+) = 1.$$

By Lemma 5.1, there exists a $\tau > 0$ such that

$$\mu^{+,\beta,(-,\varepsilon')} \preceq \mu_{\inf,\tau}^{+,\beta}$$

and therefore,

$$\mu_{\inf,\tau}^{+,\beta}(\mathcal{C}^+) = 1.$$

Therefore,

$$\Psi^{+,\beta}(\mathcal{C}_t^+ \text{ occurs for every } t \in [0,\tau]) = 1.$$

Finally, using countable additivity,

$$\Psi^{+,\beta}(\mathcal{C}_t^+ \text{ occurs for every } t) = 1. \qquad \square$$

**10. Proof of Theorem 1.4.** The aim of this section is to prove Theorem 1.4. For that we will use Theorem 1.5 and Lemma 10.1. We will not prove Lemma 10.1 since it follows immediately from the proof of Lemma 11.12 in [10] due to Y. Zhang.

A probability measure $\mu$ on $\{-1,1\}^S$ is said to have the finite energy property if all conditional probabilities on finite sets are strictly positive.

LEMMA 10.1. *Take $\mu$ to be any probability measure on $\{-1,1\}^{\mathbb{Z}^2}$ which has positive correlations and the finite energy property. Assume further that $\mu$ is invariant under translations, rotations and reflections in the coordinate axes. If $\mu(\mathcal{C}^+) = 1$, then $\mu(\mathcal{C}^-) = 0$.*

PROOF OF THEOREM 1.4. Fix $\beta > \beta_c$. By (43), there exists $\varepsilon > 0$ such that

$$\mu^{+,\beta,(-,\varepsilon)}(\mathcal{C}^+) = 1.$$

Since $\mu^{+,\beta}$ and $\pi_{1-\varepsilon}$ both have positive correlations, it follows that $\mu^{+,\beta,(-,\varepsilon)}$ has positive correlations. This is because (see [19], page 78) the product of two probability measures which have positive correlations also has positive correlations. Furthermore, a collection of increasing functions of random variables which have positive correlations also has positive correlations. In addition, the finite energy property is easily seen to hold for $\mu^{+,\beta,(-,\varepsilon)}$. Using this, we can by Lemma 10.1 conclude that

$$\mu^{+,\beta,(-,\varepsilon)}(\mathcal{C}^-) = 0.$$

By Lemma 5.1, there exists a $\tau > 0$ such that $\mu^{+,\beta,(-,\varepsilon)} \preceq \mu_{\inf,\tau}^{+,\beta}$ and hence,

$$\mu_{\inf,\tau}^{+,\beta}(\mathcal{C}^-) = 0.$$

It follows that

$$\Psi^{+,\beta}(\exists t \in [0,\tau] : \mathcal{C}_t^- \text{ occurs}) = 0,$$

and by countable additivity, we conclude

$$\Psi^{+,\beta}(\exists t \geq 0 : \mathcal{C}_t^- \text{ occurs}) = 0. \qquad \square$$

**Acknowledgment.** We thank the referee for a very careful reading and for providing a number of suggestions.

DEPARTMENT OF MATHEMATICS
CHALMERS UNIVERSITY OF TECHNOLOGY
412 96 GOTHENBURG
SWEDEN
E-MAIL: broman@math.chalmers.se
steif@math.chalmers.se
URL: www.math.chalmers.se/~broman
www.math.chalmers.se/~steif